\documentclass{article}
\usepackage{amsmath,amsthm,amsfonts,amssymb}
\usepackage{amsmath}
\usepackage{amssymb}
\usepackage{cite}
\begin{document}
\def\e#1\e{\begin{equation}#1\end{equation}}
\def\ea#1\ea{\begin{align}#1\end{align}}
\def\eq#1{{\rm(\ref{#1})}}
\theoremstyle{plain}
\newtheorem{thm}{Theorem}[section]
\newtheorem{lem}[thm]{Lemma}
\newtheorem{prop}[thm]{Proposition}
\newtheorem{conj}[thm]{Conjecture}
\newtheorem{cor}[thm]{Corollary}
\theoremstyle{definition}
\newtheorem{dfn}[thm]{Definition}
\newtheorem{prob}[thm]{Problem}
\newtheorem{ex}[thm]{Example}
\def\dim{\mathop{\rm dim}}
\def\Re{\mathop{\rm Re}}
\def\Im{\mathop{\rm Im}}
\def\ind{\mathop{\rm ind}}
\def\vol{\mathop{\rm vol}}
\def\sind{{\ts\mathop{\text{\rm s-ind}}}}
\def\supp{\mathop{\rm supp}}
\def\id{\mathop{\rm id}}
\def\U{\mathbin{\rm U}}
\def\SU{\mathop{\rm SU}}
\def\ge{\geqslant} 
\def\le{\leqslant} 
\def\N{\mathbin{\mathbb N}}
\def\R{\mathbin{\mathbb R}}
\def\Z{\mathbin{\mathbb Z}}
\def\C{\mathbin{\mathbb C}}
\def\al{\alpha}
\def\be{\beta}
\def\de{\delta}
\def\ep{\epsilon}
\def\io{\iota}
\def\la{\lambda}
\def\th{\theta}
\def\ze{\zeta}
\def\up{\upsilon}
\def\vp{\varphi}
\def\si{\sigma}
\def\om{\omega}
\def\De{\Delta}
\def\La{\Lambda}
\def\Om{\Omega}
\def\Ga{\Gamma}
\def\Si{\Sigma}
\def\Th{\Theta}
\def\Up{\Upsilon}
\def\d{{\rm d}}
\def\pd{\partial}
\def\ts{\textstyle}
\def\sst{\scriptscriptstyle}
\def\w{\wedge}
\def\sm{\setminus}
\def\op{\oplus}
\def\iy{\infty}
\def\ra{\rightarrow}
\def\longra{\longrightarrow}
\def\t{\times}
\def\na{\nabla}
\def\ha{{\textstyle\frac{1}{2}}}
\def\ti{\tilde}
\def\bs{\boldsymbol}
\def\ov{\overline}
\def\sL{{\smash{\sst L}}}
\def\sSi{{\smash{\sst\Si}}}
\def\sSii{{\smash{\sst\Si_i}}}
\def\sN{{\smash{\sst N}}}
\def\sX{{\smash{\sst X}}}
\def\sXp{{\smash{\sst X'}}}
\def\D{{\mathcal D}}
\def\I{{\mathcal I}}
\def\M{{\mathcal M}}
\def\oM{{\,\,\ov{\!\!{\mathcal M}\!}\,}}
\def\O{{\mathcal O}}
\def\cal{\mathcal}
\def\ms#1{\vert#1\vert^2}
\def\md#1{\vert #1 \vert}
\def\bmd#1{\big\vert #1 \big\vert}
\def\an#1{\langle#1\rangle}
\def\ban#1{\bigl\langle#1\bigr\rangle}
\title{Singularities of special Lagrangian submanifolds}
\author{Dominic Joyce \\ Lincoln College, Oxford}
\date{}
\maketitle

\begin{abstract}
We survey what is known about singularities of special
Lagrangian submanifolds (SL $m$-folds) in (almost)
Calabi--Yau manifolds. The bulk of the paper summarizes the
author's work \cite{Joyc13,Joyc14,Joyc15,Joyc16,Joyc17} on
SL $m$-folds $X$ with {\it isolated conical singularities}.
That is, near each singular point $x$, $X$ is modelled on
an SL cone $C$ in $\C^m$ with isolated singularity at 0.
We also discuss directions for future research, and give
a list of open problems.
\end{abstract}

\section{Introduction}
\label{gr1}

{\it Special Lagrangian $m$-folds (SL\/ $m$-folds)} are a
distinguished class of real $m$-dimensional minimal submanifolds
which may be defined in $\C^m$, or in {\it Calabi--Yau $m$-folds},
or more generally in {\it almost Calabi--Yau $m$-folds} (compact
K\"ahler $m$-folds with trivial canonical bundle). They are of
interest to Differential Geometers, to String Theorists (a
species of theoretical physicist), and perhaps in the future
to Algebraic Geometers.

This article will discuss the {\it singularities} of SL $m$-folds,
a field which has received little attention until quite recently.
We begin in \S\ref{gr2} with a brief introduction to special
Lagrangian geometry and (almost) Calabi--Yau $m$-folds.
Sections \ref{gr3}--\ref{gr7} survey the author's series
of papers \cite{Joyc13,Joyc14,Joyc15,Joyc16,Joyc17} on SL
$m$-folds with {\it isolated conical singularities}, a
large class of singularities which are simple enough to
study in detail. The last and longest section, \S\ref{gr8},
suggests directions for future research and gives some
open problems.

We say that a compact SL $m$-fold $X$ in an almost Calabi--Yau
$m$-fold $M$ for $m>2$ has {\it isolated conical singularities}
if it has only finitely many singular points $x_1,\ldots,x_n$
in $M$, such that for some {\it special Lagrangian cones} $C_i$
in $T_{\smash{x_i}}M\cong\C^m$ with $C_i\sm\{0\}$ nonsingular,
$X$ approaches $C_i$ near $x_i$, in an asymptotic $C^1$ sense.
The exact definition is given in~\S\ref{gr33}.

Section \ref{gr4} discusses the {\it regularity} of SL $m$-folds
$X$ with conical singularities $x_1,\ldots,x_n$, that is, how
quickly $X$ converges to the cone $C_i$ near $x_i$, with all
derivatives. In \S\ref{gr5} we consider the {\it deformation
theory} of compact SL $m$-folds $X$ with conical singularities.
We find that the {\it moduli space} ${\mathcal M}_\sX$ of
deformations of $X$ in $M$ is locally homeomorphic to the zeroes
of a smooth map $\Phi:{\mathcal I}_\sXp\ra{\mathcal O}_\sXp$
between finite-dimensional vector spaces, and if the {\it
obstruction space} ${\mathcal O}_\sXp$ is zero then
${\mathcal M}_\sX$ is a smooth manifold.

Section \ref{gr6} is an aside on {\it Asymptotically Conical
SL\/ $m$-folds (AC SL\/ $m$-folds)} in $\C^m$, that is,
nonsingular, noncompact SL $m$-folds $L$ in $\C^m$ which are
asymptotic at infinity to an SL cone $C$ at a prescribed
{\it rate} $\la$. In \S\ref{gr7} we explain how to {\it
desingularize} of a compact SL $m$-fold $X$ with conical
singularities $x_i$ with cones $C_i$ for $i=1,\ldots,n$ in
an almost Calabi--Yau $m$-fold $M$. We take AC SL $m$-folds
$L_i$ in $\C^m$ asymptotic to $C_i$ at infinity, and glue
$tL_i$ into $X$ at $x_i$ for small $t>0$ to get a smooth family
of compact, {\it nonsingular} SL $m$-folds $\smash{\ti N^t}$ in
$M$, with $\smash{\ti N^t}\ra X$ as~$t\ra 0$.

For brevity I generally give only statements of results, with
at most brief sketches of proofs. For the same reason I have
left out several subjects I would like to discuss. Some
particular omissions are:
\begin{itemize}
\setlength{\itemsep}{0pt}
\setlength{\parsep}{0pt}
\item We give very few {\it examples\/} of SL $m$-folds.
But many examples are known in $\C^m$, in \cite{Harv,Hask1,
HaLa,Joyc2,Joyc3,Joyc4,Joyc6,Joyc7,Joyc8,Joyc9,Joyc10,
Joyc11,Joyc12} and other papers.
\item We give no {\it applications\/} of the results of
\S\ref{gr3}--\S\ref{gr7}. See~\cite[\S 8--\S 10]{Joyc17}.
\item We do not discuss {\it smooth families} of almost
Calabi--Yau $m$-folds. However, all the main results of
\S\ref{gr24}, \S\ref{gr5} and \S\ref{gr7} have extensions
to families, which can be found in \cite{Joyc13,Joyc14,
Joyc15,Joyc16,Joyc17}. The discussion of {\it index} of
singularities in \S\ref{gr81}, and its applications in
\S\ref{gr83} and \S\ref{gr84}, would also be improved
by extending it to families.
\end{itemize}
\medskip

\noindent{\it Acknowledgements.} Many people have helped me develop
my ideas on special Lagrangian geometry. Amongst them I would like
to thank Tom Bridgeland, Adrian Butscher, Mark Haskins, Nigel Hitchin,
Stephen Marshall, Ian McIntosh, Sema Salur and Richard Thomas. I was
supported by an EPSRC Advanced Research Fellowship whilst writing
this paper.

\section{Special Lagrangian geometry}
\label{gr2}

We begin with some background from symplectic geometry.
Then special Lagrangian submanifolds (SL $m$-folds) are
introduced both in $\C^m$ and in {\it almost Calabi--Yau
$m$-folds}. We also describe the {\it deformation theory}
of compact SL $m$-folds. Some references for this section
are McDuff and Salamon \cite{McSa}, Harvey and Lawson
\cite{HaLa}, McLean \cite{McLe}, and the author~\cite{Joyc8}.

\subsection{Background from symplectic geometry}
\label{gr21}

We start by recalling some elementary symplectic geometry, which
can be found in McDuff and Salamon \cite{McSa}. Here are the basic
definitions.

\begin{dfn} Let $M$ be a smooth manifold of even dimension $2m$.
A closed $2$-form $\om$ on $M$ is called a {\it symplectic form}
if the $2m$-form $\om^m$ is nonzero at every point of $M$. Then
$(M,\om)$ is called a {\it symplectic manifold}. A submanifold
$N$ in $M$ is called {\it Lagrangian} if $\dim N=m=\ha\dim M$
and~$\om\vert_N\equiv 0$.
\label{gr2def1}
\end{dfn}

The simplest example of a symplectic manifold is~$\R^{2m}$.

\begin{dfn} Let $\R^{2m}$ have coordinates $(x_1,\ldots,x_m,
y_1,\ldots,y_m)$, and define the standard metric $g'$ and
symplectic form $\om'$ on $\R^{2m}$ by
\begin{equation*}
g'=\ts\sum_{j=1}^m(\d x_j^2+\d y_j^2)\quad\text{and}\quad
\om'=\ts\sum_{j=1}^m\d x_j\w\d y_j.
\end{equation*}
Then $(\R^{2m},\om')$ is a symplectic manifold. When we wish to
identify $\R^{2m}$ with $\C^m$, we take the complex coordinates
$(z_1,\ldots,z_m)$ on $\C^m$ to be $z_j=x_j+iy_j$. For $R>0$, define
$B_R$ to be the open ball of radius $R$ about 0 in~$\R^{2m}$.
\label{gr2def2}
\end{dfn}

{\it Darboux's Theorem} \cite[Th.~3.15]{McSa} says that every
symplectic manifold is locally isomorphic to $(\R^{2m},\om')$.
Our version easily follows.

\begin{thm} Let\/ $(M,\om)$ be a symplectic $2m$-manifold and\/
$x\in M$. Then there exist\/ $R>0$ and an embedding $\Up:B_R\ra
M$ with\/ $\Up(0)=x$ such that\/ $\Up^*(\om)=\om'$, where $\om'$
is the standard symplectic form on $\R^{2m}\supset B_R$. Given an
isomorphism $\up:\R^{2m}\!\ra\!T_xM$ with\/ $\up^*(\om\vert_x)\!=
\!\om'$, we can choose $\Up$ with\/~$\d\Up\vert_0\!=\!\up$.
\label{gr2thm1}
\end{thm}

Let $N$ be a real $m$-manifold. Then its tangent bundle $T^*N$ has
a {\it canonical symplectic form} $\hat\om$, defined as follows.
Let $(x_1,\ldots,x_m)$ be local coordinates on $N$. Extend them to
local coordinates $(x_1,\ldots,x_m,y_1,\ldots,y_m)$ on $T^*N$ such
that $(x_1,\ldots,y_m)$ represents the 1-form $y_1\d x_1+\cdots+y_m
\d x_m$ in $T_{(x_1,\ldots,x_m)}^*N$. Then~$\hat\om=\d x_1\w\d y_1+
\cdots+\d x_m\w\d y_m$.

Identify $N$ with the zero section in $T^*N$. Then $N$ is a
{\it Lagrangian submanifold\/} of $T^*N$. The {\it Lagrangian
Neighbourhood Theorem} \cite[Th.~3.33]{McSa} shows that any
compact Lagrangian submanifold $N$ in a symplectic manifold
looks locally like the zero section in~$T^*N$.

\begin{thm} Let\/ $(M,\om)$ be a symplectic manifold and\/
$N\subset M$ a compact Lagrangian submanifold. Then there
exists an open tubular neighbourhood\/ $U$ of the zero
section $N$ in $T^*N$, and an embedding $\Phi:U\ra M$ with\/
$\Phi\vert_N=\id:N\ra N$ and\/ $\Phi^*(\om)=\hat\om$, where
$\hat\om$ is the canonical symplectic structure on~$T^*N$.
\label{gr2thm2}
\end{thm}

We shall call $U,\Phi$ a {\it Lagrangian neighbourhood} of
$N$. Such neighbourhoods are useful for parametrizing nearby
Lagrangian submanifolds of $M$. Suppose that $\ti N$ is a
Lagrangian submanifold of $M$ which is $C^1$-close to $N$.
Then $\ti N$ lies in $\Phi(U)$, and is the image $\Phi\bigl(
\Ga(\al)\bigr)$ of the graph $\Ga(\al)$ of a unique
$C^1$-small 1-form $\al$ on~$N$.

As $\ti N$ is Lagrangian and $\Phi^*(\om)=\hat\om$ we see
that $\hat\om\vert_{\Ga(\al)}\equiv 0$. But one can easily
show that $\hat\om\vert_{\Ga(\al)}=-\pi^*(\d\al)$, where
$\pi:\Ga(\al)\ra N$ is the natural projection. Hence $\d\al=0$,
and $\al$ is a {\it closed\/ $1$-form}. This establishes a
1-1 correspondence between small closed 1-forms on $N$ and
Lagrangian submanifolds $\ti N$ close to $N$ in $M$, which
is an essential tool in proving later results.

\subsection{Special Lagrangian submanifolds in $\C^m$}
\label{gr22}

We define {\it calibrations} and {\it calibrated submanifolds},
following~\cite{HaLa}.

\begin{dfn} Let $(M,g)$ be a Riemannian manifold. An {\it oriented
tangent $k$-plane} $V$ on $M$ is a vector subspace $V$ of
some tangent space $T_xM$ to $M$ with $\dim V=k$, equipped
with an orientation. If $V$ is an oriented tangent $k$-plane
on $M$ then $g\vert_V$ is a Euclidean metric on $V$, so 
combining $g\vert_V$ with the orientation on $V$ gives a 
natural {\it volume form} $\vol_V$ on $V$, which is a 
$k$-form on~$V$.

Now let $\vp$ be a closed $k$-form on $M$. We say that
$\vp$ is a {\it calibration} on $M$ if for every oriented
$k$-plane $V$ on $M$ we have $\vp\vert_V\le \vol_V$. Here
$\vp\vert_V=\al\cdot\vol_V$ for some $\al\in\R$, and 
$\vp\vert_V\le\vol_V$ if $\al\le 1$. Let $N$ be an 
oriented submanifold of $M$ with dimension $k$. Then 
each tangent space $T_xN$ for $x\in N$ is an oriented
tangent $k$-plane. We say that $N$ is a {\it calibrated 
submanifold\/} if $\vp\vert_{T_xN}=\vol_{T_xN}$ for all~$x\in N$.
\label{gr2def3}
\end{dfn}

It is easy to show that calibrated submanifolds are automatically
{\it minimal submanifolds} \cite[Th.~II.4.2]{HaLa}. Here is the 
definition of special Lagrangian submanifolds in $\C^m$, taken
from~\cite[\S III]{HaLa}.

\begin{dfn} Let $\C^m$ have complex coordinates $(z_1,\dots,z_m)$, 
and define a metric $g'$, a real 2-form $\om'$ and a complex $m$-form 
$\Om'$ on $\C^m$ by
\e
\begin{split}
g'=\ms{\d z_1}+\cdots+\ms{\d z_m},\quad
\om'&=\ts\frac{i}{2}(\d z_1\w\d\bar z_1+\cdots+\d z_m\w\d\bar z_m),\\
\text{and}\quad\Om'&=\d z_1\w\cdots\w\d z_m.
\end{split}
\label{gr2eq1}
\e
Then $g',\om'$ are as in Definition \ref{gr2def2}, and $\Re\Om'$ and
$\Im\Om'$ are real $m$-forms on $\C^m$. Let $L$ be an oriented real
submanifold of $\C^m$ of real dimension $m$. We say that $L$ is a
{\it special Lagrangian submanifold\/} of $\C^m$, or {\it SL\/
$m$-fold}\/ for short, if $L$ is calibrated with respect to $
\Re\Om'$, in the sense of Definition~\ref{gr2def3}.
\label{gr2def4}
\end{dfn}

Harvey and Lawson \cite[Cor.~III.1.11]{HaLa} give the following
alternative characterization of special Lagrangian submanifolds:

\begin{prop} Let\/ $L$ be a real $m$-dimensional submanifold 
of\/ $\C^m$. Then $L$ admits an orientation making it into an
SL submanifold of\/ $\C^m$ if and only if\/ $\om'\vert_L\equiv 0$ 
and\/~$\Im\Om'\vert_L\equiv 0$.
\label{gr2prop1}
\end{prop}

Thus special Lagrangian submanifolds are {\it Lagrangian}
submanifolds satisfying the extra condition that
$\Im\Om'\vert_L\equiv 0$, which is how they get their name.

\subsection{Almost Calabi--Yau $m$-folds and SL $m$-folds} 
\label{gr23}

We shall define special Lagrangian submanifolds not just in
Calabi--Yau manifolds, as usual, but in the much larger
class of {\it almost Calabi--Yau manifolds}.

\begin{dfn} Let $m\ge 2$. An {\it almost Calabi--Yau $m$-fold\/}
is a quadruple $(M,J,\om,\Om)$ such that $(M,J)$ is a compact
$m$-dimensional complex manifold, $\om$ is the K\"ahler form
of a K\"ahler metric $g$ on $M$, and $\Om$ is a non-vanishing
holomorphic $(m,0)$-form on~$M$.

We call $(M,J,\om,\Om)$ a {\it Calabi--Yau $m$-fold\/} if in
addition $\om$ and $\Om$ satisfy
\e
\om^m/m!=(-1)^{m(m-1)/2}(i/2)^m\Om\w\bar\Om.
\label{gr2eq2}
\e
Then for each $x\in M$ there exists an isomorphism $T_xM\cong\C^m$
that identifies $g_x,\om_x$ and $\Om_x$ with the flat versions
$g',\om',\Om'$ on $\C^m$ in \eq{gr2eq1}. Furthermore, $g$ is
Ricci-flat and its holonomy group is a subgroup of~$\SU(m)$.
\label{gr2def5}
\end{dfn}

This is not the usual definition of a Calabi--Yau manifold, but
is essentially equivalent to it.

\begin{dfn} Let $(M,J,\om,\Om)$ be an almost Calabi--Yau $m$-fold,
and $N$ a real $m$-dimensional submanifold of $M$. We call $N$ a
{\it special Lagrangian submanifold}, or {\it SL $m$-fold\/} for
short, if $\om\vert_N\equiv\Im\Om\vert_N\equiv 0$. It easily
follows that $\Re\Om\vert_N$ is a nonvanishing $m$-form on $N$.
Thus $N$ is orientable, with a unique orientation in which
$\Re\Om\vert_N$ is positive.
\label{gr2def6}
\end{dfn}

Again, this is not the usual definition of SL $m$-fold, but is
essentially equivalent to it. Suppose $(M,J,\om,\Om)$ is an
almost Calabi--Yau $m$-fold, with metric $g$. Let
$\psi:M\ra(0,\iy)$ be the unique smooth function such that
\e
\psi^{2m}\om^m/m!=(-1)^{m(m-1)/2}(i/2)^m\Om\w\bar\Om,
\label{gr2eq3}
\e
and define $\ti g$ to be the conformally equivalent metric $\psi^2g$
on $M$. Then $\Re\Om$ is a {\it calibration} on the Riemannian manifold
$(M,\ti g)$, and SL $m$-folds $N$ in $(M,J,\om,\Om)$ are calibrated
with respect to it, so that they are minimal with respect to~$\ti g$.

If $M$ is a Calabi--Yau $m$-fold then $\psi\equiv 1$ by \eq{gr2eq2},
so $\ti g=g$, and an $m$-submanifold $N$ in $M$ is special Lagrangian
if and only if it is calibrated w.r.t.\ $\Re\Om$ on $(M,g)$, as in
Definition \ref{gr2def4}. This recovers the usual definition of
special Lagrangian $m$-folds in Calabi--Yau $m$-folds.

\subsection{Deformations of compact SL $m$-folds} 
\label{gr24}

The {\it deformation theory} of special Lagrangian submanifolds
was studied by McLean \cite[\S 3]{McLe}, who proved the following
result in the Calabi--Yau case. The extension to the almost
Calabi--Yau case is described in~\cite[\S 9.5]{Joyc8}.

\begin{thm} Let\/ $N$ be a compact SL\/ $m$-fold in an almost
Calabi--Yau $m$-fold\/ $(M,J,\om,\Om)$. Then the moduli space
$\M_\sN$ of special Lagrangian deformations of\/ $N$ is a smooth
manifold of dimension~$b^1(N)$.
\label{gr2thm3}
\end{thm}

\begin{proof}[Sketch of proof] There is a natural orthogonal
decomposition $TM\vert_N=TN\op\nu$, where $\nu\ra N$ is the
{\it normal bundle} of $N$ in $M$. As $N$ is Lagrangian, the
complex structure $J:TM\ra TM$ gives an isomorphism $J:\nu\ra TN$.
But the metric $g$ gives an isomorphism $TN\cong T^*N$. Composing
these two gives an isomorphism~$\nu\cong T^*N$.

Let $T$ be a small {\it tubular neighborhood} of $N$ in $M$. Then 
we can identify $T$ with a neighborhood of the zero section in $\nu$.
Using the isomorphism $\nu\cong T^*N$, we have an identification
between $T$ and a neighborhood of the zero section in $T^*N$. This
can be chosen to identify the K\"ahler form $\om$ on $T$ with the natural
symplectic structure on $T^*N$. Let $\pi:T\ra N$ be the obvious projection.

Under this identification, submanifolds $N'$ in $T\subset M$ which 
are $C^1$ close to $N$ are identified with the graphs of small smooth 
sections $\al$ of $T^*N$. That is, submanifolds $N'$ of $M$ close to
$N$ are identified with 1-{\it forms} $\al$ on $N$. We need to know: 
which 1-forms $\al$ are identified with {\it special Lagrangian}
submanifolds~$N'$?

Well, $N'$ is special Lagrangian if $\om\vert_{N'}\equiv
\Im\Om\vert_{N'}\equiv 0$. Now $\pi\vert_{N'}:N'\ra N$ is a
diffeomorphism, so we can push $\om\vert_{N'}$ and
$\Im\Om\vert_{N'}$ down to $N$, and regard them as functions 
of $\al$. Calculation shows that $\pi_*\bigl(\om\vert_{N'}\bigr)
=\d\al$ and $\pi_*\bigl(\Im\Om\vert_{N'}\bigr)=F(\al,\nabla\al)$,
where $F$ is a nonlinear function of its arguments. Thus, the moduli 
space ${\cal M}_\sN$ is locally isomorphic to the set of small 1-forms 
$\al$ on $N$ such that $\d\al\equiv 0$ and $F(\al,\nabla\al)\equiv 0$.

Now it turns out that $F$ satisfies $F(\al,\nabla\al)\approx \d(*\al)$ 
when $\al$ is small. Therefore ${\cal M}_\sN$ is locally approximately 
isomorphic to the vector space of 1-forms $\al$ with $\d\al=\d(*\al)=0$.
But by Hodge theory, this is isomorphic to the de Rham cohomology
group $H^1(N,\R)$, and is a manifold with dimension~$b^1(N)$.

To carry out this last step rigorously requires some technical
machinery: one must work with certain {\it Banach spaces} of 
sections of $T^*N$, $\La^2T^*N$ and $\La^mT^*N$, use {\it elliptic 
regularity results} to show the map $\al\mapsto\bigl(\d\al,
F(\al,\nabla\al)\bigr)$ has {\it closed image} in these Banach spaces,
and then use the {\it Implicit Function Theorem for Banach spaces}
to show that the kernel of the map is what we expect.
\end{proof}

\section{SL cones and conical singularities}
\label{gr3}

We begin in \S\ref{gr31} with some definitions on {\it special
Lagrangian cones}. Section \ref{gr32} gives {\it examples} of
SL cones, and \S\ref{gr33} defines {\it SL\/ $m$-folds with
conical singularities}, the subject of the paper. Section
\ref{gr34} discusses {\it homology} and {\it cohomology} of
SL $m$-folds with conical singularities.

\subsection{Preliminaries on special Lagrangian cones}
\label{gr31}

We define {\it special Lagrangian cones}, and some notation.

\begin{dfn} A (singular) SL $m$-fold $C$ in $\C^m$ is called a
{\it cone} if $C=tC$ for all $t>0$, where $tC=\{t\,{\bf x}:{\bf x}
\in C\}$. Let $C$ be a closed SL cone in $\C^m$ with an isolated
singularity at 0. Then $\Si=C\cap{\mathcal S}^{2m-1}$ is a compact,
nonsingular $(m\!-\!1)$-submanifold of ${\mathcal S}^{2m-1}$, not
necessarily connected. Let $g_\sSi$ be the restriction
of $g'$ to $\Si$, where $g'$ is as in~\eq{gr2eq1}.

Set $C'=C\sm\{0\}$. Define $\io:\Si\t(0,\iy)\ra\C^m$ by
$\io(\si,r)=r\si$. Then $\io$ has image $C'$. By an abuse
of notation, {\it identify} $C'$ with $\Si\t(0,\iy)$ using
$\io$. The {\it cone metric} on $C'\cong\Si\t(0,\iy)$
is~$g'=\io^*(g')=\d r^2+r^2g_\sSi$.

For $\al\in\R$, we say that a function $u:C'\ra\R$ is
{\it homogeneous of order} $\al$ if $u\circ t\equiv t^\al u$ for
all $t>0$. Equivalently, $u$ is homogeneous of order $\al$ if
$u(\si,r)\equiv r^\al v(\si)$ for some function~$v:\Si\ra\R$.
\label{gr3def1}
\end{dfn}

In \cite[Lem.~2.3]{Joyc13} we study {\it homogeneous harmonic
functions} on~$C'$.

\begin{lem} In the situation of Definition \ref{gr3def1},
let\/ $u(\si,r)\equiv r^\al v(\si)$ be a homogeneous function
of order $\al$ on $C'=\Si\t(0,\iy)$, for $v\in C^2(\Si)$. Then
\begin{equation*}
\De u(\si,r)=r^{\al-2}\bigl(\De_\sSi v-\al(\al+m-2)v\bigr),
\end{equation*}
where $\De$, $\De_\sSi$ are the Laplacians on $(C',g')$
and\/ $(\Si,g_\sSi)$. Hence, $u$ is harmonic on $C'$
if and only if\/ $v$ is an eigenfunction of\/ $\De_\sSi$
with eigenvalue~$\al(\al+m-2)$.
\label{gr3lem}
\end{lem}

Following \cite[Def.~2.5]{Joyc13}, we define:

\begin{dfn} In Definition \ref{gr3def1}, suppose $m>2$ and define
\e
\D_\sSi=\bigl\{\al\in\R:\text{$\al(\al+m-2)$ is
an eigenvalue of $\De_\sSi$}\bigr\}.
\label{gr3eq1}
\e
Then $\D_\sSi$ is a countable, discrete subset of
$\R$. By Lemma \ref{gr3lem}, an equivalent definition is that
$\D_\sSi$ is the set of $\al\in\R$ for which there
exists a nonzero homogeneous harmonic function $u$ of order
$\al$ on~$C'$.

Define $m_\sSi:\D_\sSi\ra\N$ by taking
$m_\sSi(\al)$ to be the multiplicity of the eigenvalue
$\al(\al+m-2)$ of $\De_\sSi$, or equivalently the
dimension of the vector space of homogeneous harmonic
functions $u$ of order $\al$ on $C'$. Define
$N_\sSi:\R\ra\Z$ by
\e
N_\sSi(\de)=
-\sum_{\!\!\!\!\al\in\D_\sSi\cap(\de,0)\!\!\!\!}m_\sSi(\al)
\;\>\text{if $\de<0$, and}\;\>
N_\sSi(\de)=
\sum_{\!\!\!\!\al\in\D_\sSi\cap[0,\de]\!\!\!\!}m_\sSi(\al)
\;\>\text{if $\de\ge 0$.}
\label{gr3eq2}
\e
Then $N_\sSi$ is monotone increasing and upper semicontinuous,
and is discontinuous exactly on $\D_\sSi$, increasing by
$m_\sSi(\al)$ at each $\al\in\D_\sSi$. As the
eigenvalues of $\De_\sSi$ are nonnegative, we see that
$\D_\sSi\cap(2-m,0)=\emptyset$ and $N_\sSi\equiv 0$
on~$(2-m,0)$.
\label{gr3def2}
\end{dfn}

We define the {\it stability index} of $C$, and {\it stable}
and {\it rigid} cones~\cite[Def.~3.6]{Joyc14}.

\begin{dfn} Let $C$ be an SL cone in $\C^m$ for $m>2$ with an
isolated singularity at 0, let $G$ be the Lie subgroup of\/
$\SU(m)$ preserving $C$, and use the notation of Definitions
\ref{gr3def1} and \ref{gr3def2}. Then \cite[eq.~(8)]{Joyc14}
shows that
\e
m_\sSi(0)=b^0(\Si),\quad
m_\sSi(1)\ge 2m \quad\text{and}\quad
m_\sSi(2)\ge m^2-1-\dim G.
\label{gr3eq3}
\e

Define the {\it stability index} $\sind(C)$ to be
\begin{equation*}
\sind(C)=N_\sSi(2)-b^0(\Si)-m^2-2m+1+\dim G.
\end{equation*}
Then $\sind(C)\ge 0$ by \eq{gr3eq3}, as $N_\sSi(2)\ge
m_\sSi(0)+m_\sSi(1)+m_\sSi(2)$ by \eq{gr3eq2}.
We call $C$ {\it stable} if~$\sind(C)=0$.

Following \cite[Def.~6.7]{Joyc13}, we call $C$ {\it rigid} if
$m_\sSi(2)=m^2-1-\dim G$. As
\begin{equation*}
\sind(C)\ge m_\sSi(2)-(m^2-1-\dim G)\ge 0,
\end{equation*}
we see that {\it if\/ $C$ is stable, then $C$ is rigid}.
\label{gr3def3}
\end{dfn}

We shall see in \S\ref{gr5} that $\sind(C)$ is the dimension
of an obstruction space to deforming an SL $m$-fold $X$ with
a conical singularity with cone $C$, and that if $C$ is
{\it stable} then the deformation theory of $X$ simplifies.
An SL cone $C$ is {\it rigid} if all infinitesimal deformations
of $C$ as an SL cone come from $\SU(m)$ rotations of $C$. This
will be useful in the Geometric Measure Theory material
of~\S\ref{gr4}.

\subsection{Examples of special Lagrangian cones}
\label{gr32}

In our first example we can compute the data of \S\ref{gr31}
very explicitly.

\begin{ex} Here is a family of special Lagrangian cones
constructed by Harvey and Lawson \cite[\S III.3.A]{HaLa}.
For $m\ge 3$, define
\e
C_{\sst\rm HL}^m=\bigl\{(z_1,\ldots,z_m)\in\C^m:i^{m+1}z_1\cdots
z_m\in[0,\iy), \quad \md{z_1}=\cdots=\md{z_m}\bigr\}.
\label{gr3eq4}
\e
Then $C_{\sst\rm HL}^m$ is a {\it special Lagrangian cone} in
$\C^m$ with an isolated singularity at 0, and $\Si_{\sst\rm HL}^m
=C_{\sst\rm HL}^m\cap{\mathcal S}^{2m-1}$ is an $(m\!-\!1)$-torus
$T^{m-1}$. Both $C_{\sst\rm HL}^m$ and $\Si_{\sst\rm HL}^m$ are
invariant under the $\U(1)^{m-1}$ subgroup of $\SU(m)$ acting by
\e
(z_1,\ldots,z_m)\mapsto({\rm e}^{i\th_1}z_1,\ldots,{\rm e}^{i\th_m}
z_m) \quad\text{for $\th_j\in\R$ with $\th_1+\cdots+\th_m=0$.}
\label{gr3eq5}
\e
In fact $\pm\,C_{\sst\rm HL}^m$ for $m$ odd, and $C_{\sst\rm
HL}^m,iC_{\sst\rm HL}^m$ for $m$ even, are the unique SL cones in
$\C^m$ invariant under \eq{gr3eq5}, which is how Harvey and Lawson
constructed them.

The metric on $\Si_{\sst\rm HL}^m\cong T^{m-1}$ is flat, so it is
easy to compute the eigenvalues of $\De_{\smash{\sst\Si_{\rm HL}^m}}$.
This was done by Marshall \cite[\S 6.3.4]{Mars}. There is a 1-1
correspondence between $(n_1,\ldots,n_{m-1})\in\Z^{m-1}$ and
eigenvectors of $\De_{\smash{\sst\Si_{\rm HL}^m}}$ with eigenvalue
\e
m\sum_{i=1}^{m-1}n_i^2-\sum_{i,j=1}^{m-1}n_in_j.
\label{gr3eq6}
\e

Using \eq{gr3eq6} and a computer we can find the eigenvalues
of $\De_{\smash{\sst\Si_{\rm HL}^m}}$ and their multiplicities.
The Lie subgroup $G_{\sst\rm HL}^m$ of $\SU(m)$ preserving
$C_{\sst\rm HL}^m$ has identity component the $\U(1)^{m-1}$ of
\eq{gr3eq5}, so that $\dim G_{\sst\rm HL}^m=m-1$. Thus we can
calculate $\D_{\smash{\sst\Si_{\rm HL}^m}}$, $m_{\smash{\sst
\Si_{\rm HL}^m}}$, $N_{\smash{\sst\Si_{\rm HL}^m}}$, and the
stability index $\sind(C_{\sst\rm HL}^m)$. This was done by
Marshall \cite[Table 6.1]{Mars} and the author
\cite[\S 3.2]{Joyc14}. Table \ref{gr3table} gives the data
$m,N_{\smash{\sst\Si_{\rm HL}^m}}(2),m_{\smash{\sst\Si_{\rm
HL}^m}}(2)$ and $\sind(C_{\sst\rm HL}^m)$ for~$3\le m\le 12$.

\begin{table}[htb]
\center{
\begin{tabular}{|l|r|r|r|r|r|r|r|r|r|r|}\hline
$m\vphantom{\bigr(^l_j}$
&  3 &  4 &  5 &  6 &  7 &  8 &  9 &  10 & 11 & 12 \\
\hline
$N_{\smash{\sst\Si_{\rm HL}^m}}(2)\vphantom{\bigr(^l_j}$
& 13 & 27 & 51 & 93 & 169& 311& 331& 201 & 243& 289\\
\hline
$m_{\smash{\sst\Si_{\rm HL}^m}}(2)\vphantom{\bigr(^l_j}$
&  6 & 12 & 20 & 30 & 42 & 126& 240&  90 & 110& 132\\
\hline
$\sind(C_{\sst\rm HL}^m)\vphantom{\bigr(^l_j}$
&  0 &  6 & 20 & 50 & 112& 238& 240&  90 & 110& 132\\
\hline
\end{tabular}
}
\caption{Data for $\U(1)^{m-1}$-invariant SL cones
$C_{\sst\rm HL}^m$ in $\C^m$}
\label{gr3table}
\end{table}

One can also prove that
\e
N_{\smash{\sst\Si_{\rm HL}^m}}(2)=2m^2+1\;\>\text{and}\;\>
m_{\smash{\sst\Si_{\rm HL}^m}}(2)=\sind(C_{\sst\rm HL}^m)=m^2-m
\;\>\text{for $m\ge 10$.}
\label{gr3eq7}
\e
As $C_{\sst\rm HL}^m$ is {\it stable} when $\sind(C_{\sst\rm
HL}^m)=0$ we see from Table \ref{gr3table} and \eq{gr3eq7}
that $C_{\sst\rm HL}^3$ is a {\it stable} cone in $\C^3$, but
$C_{\sst\rm HL}^m$ is {\it unstable} for $m\ge 4$. Also
$C_{\sst\rm HL}^m$ is {\it rigid\/} when $m_{\smash{\sst\Si_{
\rm HL}^m}}(2)=m^2-m$, as $\dim G_{\sst\rm HL}^m=m-1$. Thus
$C_{\sst\rm HL}^m$ is {\it rigid\/} if and only if $m\ne 8,9$,
by Table \ref{gr3table} and~\eq{gr3eq7}.
\label{gr3ex1}
\end{ex}

Here is an example chosen from \cite[Ex.~9.4]{Joyc2} as it is
easy to write down.

\begin{ex} Let $a_1,\ldots,a_m\in\Z$ with $a_1+\cdots+a_m=0$
and highest common factor 1, such that $a_1,\ldots,a_k>0$
and $a_{k+1},\ldots,a_m<0$ for $0<k<m$. Define
\begin{align*}
L^{a_1,\ldots,a_m}_0=\bigl\{
\bigl(i{\rm e}^{ia_1\th}x_1&,{\rm e}^{ia_2\th}x_2,\ldots,
{\rm e}^{ia_m\th}x_m\bigr):
\th\in[0,2\pi),\\ 
&x_1,\ldots,x_m\in\R,\qquad 
a_1x_1^2+\cdots+a_mx_m^2=0\bigr\}.
\end{align*}
Then $L^{a_1,\ldots,a_m}_0$ is an {\it immersed SL cone} in
$\C^m$, with an isolated singularity at~0.

Define $C^{a_1,\ldots,a_m}=\bigl\{(x_1,\ldots,x_m)\in\R^m:
a_1x_1^2+\cdots+a_mx_m^2=0\bigr\}$. Then $C^{a_1,\ldots,a_m}$
is a quadric cone on ${\mathcal S}^{k-1}\t S^{m-k-1}$ in $\R^m$,
and $L^{a_1,\ldots,a_m}_0$ is the image of an immersion $\Phi:
C^{a_1,\ldots,a_m}\t{\mathcal S}^1\ra\C^m$, which is generically
2:1. Therefore $L^{a_1,\ldots,a_m}_0$ is an immersed SL cone
on~$({\mathcal S}^{k-1}\t{\mathcal S}^{m-k-1}\t{\mathcal S}^1)/\Z_2$.
\label{gr3ex2}
\end{ex}

Further examples of SL cones are constructed by Harvey and
Lawson \cite[\S III.3]{HaLa}, Haskins \cite{Hask1}, the
author \cite{Joyc2,Joyc3}, and others. Special Lagrangian
cones in $\C^3$ are a special case, which may be treated
using the theory of {\it integrable systems}. In principle
this should yield a {\it classification} of all SL cones on
$T^2$ in $\C^3$. For more information see McIntosh \cite{McIn}
or the author~\cite{Joyc7}.

\subsection{Special Lagrangian $m$-folds with conical singularities}
\label{gr33}

Now we can define {\it conical singularities} of SL $m$-folds,
following~\cite[Def.~3.6]{Joyc13}.

\begin{dfn} Let $(M,J,\om,\Om)$ be an almost Calabi--Yau $m$-fold
for $m>2$, and define $\psi:M\ra(0,\iy)$ as in \eq{gr2eq3}. Suppose
$X$ is a compact singular SL $m$-fold in $M$ with singularities at
distinct points $x_1,\ldots,x_n\in X$, and no other singularities.

Fix isomorphisms $\up_i:\C^m\ra T_{x_i}M$ for $i=1,\ldots,n$
such that $\up_i^*(\om)=\om'$ and $\up_i^*(\Om)=\psi(x_i)^m\Om'$,
where $\om',\Om'$ are as in \eq{gr2eq1}. Let $C_1,\ldots,C_n$ be SL
cones in $\C^m$ with isolated singularities at 0. For $i=1,\ldots,n$
let $\Si_i=C_i\cap{\mathcal S}^{2m-1}$, and let $\mu_i\in(2,3)$ with
\e
(2,\mu_i]\cap\D_\sSii=\emptyset,
\quad\text{where $\D_\sSii$ is defined in \eq{gr3eq1}.}
\label{gr3eq8}
\e
Then we say that
$X$ has a {\it conical singularity} or {\it conical singular
point} at $x_i$, with {\it rate} $\mu_i$ and {\it cone} $C_i$
for $i=1,\ldots,n$, if the following holds.

By Theorem \ref{gr2thm1} there exist embeddings $\Up_i:B_R\ra M$
for $i=1,\ldots,n$ satisfying $\Up_i(0)=x_i$, $\d\Up_i\vert_0=\up_i$
and $\Up_i^*(\om)=\om'$, where $B_R$ is the open ball of radius $R$
about 0 in $\C^m$ for some small $R>0$. Define $\io_i:\Si_i\t(0,R)
\ra B_R$ by $\io_i(\si,r)=r\si$ for~$i=1,\ldots,n$.

Define $X'=X\sm\{x_1,\ldots,x_n\}$. Then there should exist a
compact subset $K\subset X'$ such that $X'\sm K$ is a union of
open sets $S_1,\ldots,S_n$ with $S_i\subset\Up_i(B_R)$, whose
closures $\bar S_1,\ldots,\bar S_n$ are disjoint in $X$. For
$i=1,\ldots,n$ and some $R'\in(0,R]$ there should exist a smooth
$\phi_i:\Si_i\t(0,R')\ra B_R$ such that $\Up_i\circ\phi_i:\Si_i
\t(0,R')\ra M$ is a diffeomorphism $\Si_i\t(0,R')\ra S_i$, and
\e
\bmd{\na^k(\phi_i-\io_i)}=O(r^{\mu_i-1-k})
\quad\text{as $r\ra 0$ for $k=0,1$.}
\label{gr3eq9}
\e
Here $\na$ is the Levi-Civita connection of the cone metric
$\io_i^*(g')$ on $\Si_i\t(0,R')$, $\md{\,.\,}$ is computed
using $\io_i^*(g')$. If the cones $C_1,\ldots,C_n$ are
{\it stable} in the sense of Definition \ref{gr3def3}, then
we say that $X$ has {\it stable conical singularities}.
\label{gr3def4}
\end{dfn}

We will see in Theorems \ref{gr4thm1} and \ref{gr4thm2} that
if \eq{gr3eq9} holds for $k=0,1$ and some $\mu_i$ satisfying
\eq{gr3eq8}, then we can choose a natural $\phi_i$ for which
\eq{gr3eq9} holds for {\it all\/} $k\ge 0$, and for {\it all\/}
rates $\mu_i$ satisfying \eq{gr3eq8}. Thus the number of
derivatives required in \eq{gr3eq9} and the choice of $\mu_i$
both make little difference. We choose $k=0,1$ in \eq{gr3eq9},
and some $\mu_i$ in \eq{gr3eq8}, to make the definition as weak
as possible.

We suppose $m>2$ for two reasons. Firstly, the only SL cones
in $\C^2$ are finite unions of SL planes $\R^2$ in $\C^2$
intersecting only at 0. Thus any SL 2-fold with conical
singularities is actually {\it nonsingular} as an immersed
2-fold, so there is really no point in studying them.
Secondly, $m=2$ is a special case in the analysis of
\cite[\S 2]{Joyc13}, and it is simpler to exclude it.
Therefore we will suppose $m>2$ throughout the paper.

Here are the reasons for the conditions on $\mu_i$ in
Definition~\ref{gr3def4}:
\begin{itemize}
\setlength{\itemsep}{0pt}
\setlength{\parsep}{0pt}
\item We need $\mu_i>2$, or else \eq{gr3eq9} does not force
$X$ to approach $C_i$ near~$x_i$.
\item The definition involves a choice of $\Up_i:B_R\ra M$.
If we replace $\Up_i$ by a different choice $\ti\Up_i$ then
we should replace $\phi_i$ by $\ti\phi_i=(\ti\Up_i^{-1}\circ
\Up_i)\circ\phi_i$ near 0 in $B_R$. Calculation shows that
as $\Up_i,\ti\Up_i$ agree up to second order, we
have~$\bmd{\na^k(\ti\phi_i-\phi_i)}=O(r^{2-k})$.

Therefore we choose $\mu_i<3$ so that these $O(r^{2-k})$
terms are absorbed into the $O(r^{\mu_i-1-k})$ in \eq{gr3eq9}.
This makes the definition independent of the choice of
$\Up_i$, which it would not be if~$\mu_i>3$.

\item Condition \eq{gr3eq8} is needed to prove the regularity
result Theorem \ref{gr4thm2}, and also to reduce to a minimum
the obstructions to deforming compact SL $m$-folds with
conical singularities studied in~\S\ref{gr5}.
\end{itemize}

\subsection{Homology and cohomology}
\label{gr34}

Next we discuss {\it homology} and {\it cohomology} of SL $m$-folds
with conical singularities, following \cite[\S 2.4]{Joyc13}. For a
general reference, see for instance Bredon \cite{Bred}. When $Y$
is a manifold, write $H^k(Y,\R)$ for the $k^{\rm th}$ {\it de Rham
cohomology group} and $H^k_{\rm cs}(Y,\R)$ for the $k^{\rm th}$
{\it compactly-supported de Rham cohomology group} of $Y$. If $Y$
is compact then $H^k(Y,\R)=H^k_{\rm cs}(Y,\R)$. The {\it Betti
numbers} of $Y$ are $b^k(Y)=\dim H^k(Y,\R)$ and~$b^k_{\rm cs}(Y)
=\dim H^k_{\rm cs}(Y,\R)$.

Let $Y$ be a topological space, and $Z\subset Y$ a subspace.
Write $H_k(Y,\R)$ for the $k^{\rm th}$ {\it real singular
homology group} of $Y$, and $H_k(Y;Z,\R)$ for the $k^{\rm th}$
{\it real singular relative homology group} of $(Y;Z)$. When
$Y$ is a manifold and $Z$ a submanifold we define $H_k(Y,\R)$
and $H_k(Y;Z,\R)$ using {\it smooth\/} simplices, as in
\cite[\S V.5]{Bred}. Then the pairing between (singular)
homology and (de Rham) cohomology is defined at the chain
level by integrating $k$-forms over $k$-simplices.

Let $X$ be a compact SL $m$-fold in $M$ with conical
singularities $x_1,\ldots,x_n$ and cones $C_1,\ldots,C_n$, and
set $X'=X\sm\{x_1,\ldots,x_n\}$ and $\Si_i=C_i\cap{\mathcal S}^{2m-1}$,
as in \S\ref{gr33}. Then $X'$ is the interior of a compact manifold
$\bar X'$ with boundary $\coprod_{i=1}^n\Si_i$. Using this we show
in \cite[\S 2.4]{Joyc13} that there is a natural long exact sequence
\e
\cdots\!\ra\!
H^k_{\rm cs}(X',\R)\!\ra\!H^k(X',\R)\!\ra\!\ts\bigoplus_{i=1}^n
H^k(\Si_i,\R)\!\ra\!H^{k+1}_{\rm cs}(X',\R)\!\ra\!\cdots,
\label{gr3eq10}
\e
and natural isomorphisms
\begin{gather*}
H_k\bigl(X;\{x_1,\ldots,x_n\},\R\bigr)^*\!\cong\!
H^k_{\rm cs}(X',\R)\!\cong\!H_{m-k}(X',\R)\!\cong\!H^{m-k}(X',\R)^*\\
\text{and}\quad
H^k_{\rm cs}(X',\R)\cong H_k(X,\R)^*
\quad\text{for all $k>1$.}
\end{gather*}
The inclusion $\io:X\ra M$ induces homomorphisms~$\io_*:
H_k(X,\R)\ra H_k(M,\R)$.

\section{The asymptotic behaviour of $X$ near $x_i$}
\label{gr4}

We now review the work of \cite{Joyc13} on the {\it regularity}
of SL $m$-folds with conical singularities. Let $M$ be an
almost Calabi--Yau $m$-fold and $X$ an SL $m$-fold in $M$
with conical singularities at $x_1,\ldots,x_n$, with
identifications $\up_i$ and cones $C_i$. We study how
quickly $X$ converges to the cone $\up(C_i)$ in
$T_{\smash{x_i}}M$ near~$x_i$.

Roughly speaking, we work by arranging for $\phi_i$ in
Definition \ref{gr3def4} to satisfy an {\it elliptic equation},
and then use {\it elliptic regularity results} to deduce
asymptotic bounds for $\phi_i-\io_i$ and all its derivatives.
Now $\phi_i$ is not uniquely defined, but is a more-or-less
arbitrary parametrization of $\Up_i^*(X')$ near 0 in $\C^m$.
To make $\phi_i$ satisfy an elliptic equation we impose an
{\it extra condition}, that $(\phi_i-\io_i)(\si,r)$ is
orthogonal to $T_{\io_i(\si,r)}C_i$ w.r.t.\ the metric $g'$
on $\C^m$, for all $(\si,r)\in\Si_i\t(0,R')$. By
\cite[Th.~4.4]{Joyc13} this also fixes $\phi_i$ uniquely,
given $\up_i,R,\Up_i$ and~$R'$.

\begin{thm} Let\/ $(M,J,\om,\Om)$ be an almost Calabi--Yau
$m$-fold, and\/ $X$ a compact SL $m$-fold in $M$ with conical
singularities at\/ $x_1,\ldots,x_n$ with identifications
$\up_i:\C^m\ra T_{x_i}M$ and cones $C_1,\ldots,C_n$.
Choose $R>0$ and\/ $\Up_i:B_R\ra M$ as in Definition
\ref{gr3def4}. Then for sufficiently small\/ $R'\in(0,R]$
there exist unique $\phi_i,S_i$ for $i=1,\ldots,n$
satisfying the conditions of Definition \ref{gr3def4} and
\e
(\phi_i-\io_i)(\si,r)\perp T_{\io_i(\si,r)}C_i
\quad\text{in $\C^m$ for all $(\si,r)\in\Si_i\t(0,R')$.}
\label{gr4eq1}
\e
\label{gr4thm1}
\end{thm}

In fact \cite[Th.~4.4]{Joyc13} characterizes $\phi_i$ in terms of a
{\it Lagrangian neighbourhood\/} $U_{\sst C_i},\Phi_{\sst C_i}$ of
$C_i$ in $\C^m$, but examining the proof of \cite[Th.~4.2]{Joyc13}
shows this is equivalent to \eq{gr4eq1}. In \cite[\S 5]{Joyc13} we
study the asymptotic behaviour of the maps $\phi_i$ of Theorem
\ref{gr4thm1}. Combining \cite[Th.s 5.1 \& 5.5, Lem.~4.5]{Joyc13}
proves:

\begin{thm} In the situation of Theorem \ref{gr4thm1}, suppose
$\mu_i'\in(2,3)$ with\/ $(2,\mu_i']\cap\D_\sSii=\emptyset$ for
$i=1,\ldots,n$. Then
\e
\begin{gathered}
\bmd{\na^k(\phi_i-\io_i)}=O(r^{\mu_i'-1-k})\quad
\quad\text{for all\/ $k\ge 0$ and\/ $i=1,\ldots,n$.}
\end{gathered}
\label{gr4eq2}
\e

Hence $X$ has conical singularities at $x_i$ with cone $C_i$
and rate $\mu_i'$, for all possible rates $\mu_i'$ allowed by
Definition \ref{gr3def4}. Therefore, the definition of
conical singularities is essentially independent of the
choice of rate~$\mu_i$.
\label{gr4thm2}
\end{thm}

Theorem \ref{gr4thm2} in effect {\it strengthens} the definition of
SL $m$-folds with conical singularities, Definition \ref{gr3def4},
as it shows that \eq{gr3eq9} actually implies the much stronger
condition \eq{gr4eq2} on all derivatives.

To prove Theorem \ref{gr4thm2}, we show using an analogue of
Theorem \ref{gr2thm2} for $C_i$ in $\C^m$ that as $\Up_i^*(X')$
is {\it Lagrangian} in $B_R$, we may regard $\phi_i$ as the
graph of a {\it closed\/ $1$-form} $\eta_i$ on $\Si_i\t(0,R')$.
The asymptotic condition \eq{gr3eq9} implies that $\eta_i$ is
{\it exact}, so we may write $\eta_i=\d A_i$ for smooth
$A_i:\Si_i\t(0,R')\ra\R$. As $\Im\Om\vert_{X'}\equiv 0$, we
find that $A_i$ satisfies the second-order nonlinear~p.d.e.
\e
\d^*\bigl(\psi^m\d A_i\bigr)(\si,r)=
Q\bigl(\si,r,\d A_i(\si,r),\na^2A_i(\si,r)\bigr)
\label{gr4eq3}
\e
for $(\si,r)\in\Si_i\t(0,R')$, where $Q$ is a smooth nonlinear
function.

When $r$ is small the $Q$ term in \eq{gr4eq3} is also small
and \eq{gr4eq3} approximates $\De_iA_i=0$, where $\De_i$ is
the Laplacian on the cone $C_i$. Therefore \eq{gr4eq3} is
{\it elliptic} for small $r$. Using known results on the
regularity of solutions of nonlinear second-order elliptic
p.d.e.s, and a theory of analysis on weighted Sobolev spaces
on manifolds with cylindrical ends developed by Lockhart and
McOwen \cite{Lock}, we can then prove~\eq{gr4eq2}.

Our next result \cite[Th.~6.8]{Joyc13} is an application of
{\it Geometric Measure Theory}. For an introduction to the
subject, see Morgan \cite{Morg}. Geometric Measure Theory
studies measure-theoretic generalizations of submanifolds
called {\it integral currents}, which may be very singular,
and is particularly powerful for {\it minimal\/} submanifolds.
As from \S\ref{gr2} SL $m$-folds are minimal submanifolds
w.r.t.\ an appropriate metric, many major results of Geometric
Measure Theory immediately apply to {\it special Lagrangian
integral currents}, a very general class of singular SL
$m$-folds with strong compactness properties.

\begin{thm} Let\/ $(M,J,\om,\Om)$ be an almost Calabi--Yau
$m$-fold and define $\psi:M\ra(0,\iy)$ as in \eq{gr2eq3}.
Let\/ $x\in M$ and fix an isomorphism $\up:\C^m\ra T_xM$
with\/ $\up^*(\om)=\om'$ and\/ $\up^*(\Om)=\psi(x)^m\Om'$,
where $\om',\Om'$ are as in~\eq{gr2eq1}.

Suppose that\/ $T$ is a special Lagrangian integral current in
$M$ with\/ $x\in T^\circ$, where $T^\circ=\supp T\sm\supp\pd T$,
and that\/ $\up_*(C)$ is a multiplicity $1$ tangent cone to $T$
at\/ $x$, where $C$ is a rigid special Lagrangian cone in $\C^m$,
in the sense of Definition \ref{gr3def3}. Then $T$ has a conical
singularity at\/ $x$, in the sense of Definition~\ref{gr3def4}.
\label{gr4thm3}
\end{thm}

This is a {\it weakening} of Definition \ref{gr3def4} for
{\it rigid\/} cones $C$. Theorem \ref{gr4thm3} also holds
for the larger class of {\it Jacobi integrable} SL cones
$C$, defined in \cite[Def.~6.7]{Joyc13}. Basically, Theorem
\ref{gr4thm3} shows that if a singular SL $m$-fold $T$ in $M$
is locally modelled on a rigid SL cone $C$ in only a very
weak sense, then it necessarily satisfies Definition
\ref{gr3def4}. One moral of Theorems \ref{gr4thm2} and
\ref{gr4thm3} is that, at least for rigid SL cones $C$,
more-or-less {\it any} sensible definition of SL $m$-folds
with conical singularities is equivalent to Definition~\ref{gr3def4}.

Theorem \ref{gr4thm3} is proved by applying regularity
results of Allard and Almgren, and Adams and Simon, mildly
adapted to the special Lagrangian situation, which roughly
say that if a tangent cone $C_i$ to $X$ at $x_i$ has an
isolated singularity at 0, is multiplicity 1, and rigid,
then $X$ has a parametrization $\phi_i$ near $x_i$ which
satisfies \eq{gr3eq9} for some $\mu_i>2$. It then quickly
follows that $X$ has a conical singularity at $x_i$, in
the sense of Definition~\ref{gr3def4}.

As discussed in \cite[\S 6.3]{Joyc13}, one can use other
results from Geometric Measure Theory to argue that for
tangent cones $C$ which are not Jacobi integrable, Definition
\ref{gr3def4} may be {\it too strong}, in that there could
exist examples of singular SL $m$-folds with tangent cone
$C$ which are not covered by Definition \ref{gr3def4}, as
the decay conditions \eq{gr3eq9} are too strict. 

\section{Moduli of SL $m$-folds with conical singularities}
\label{gr5}

Next we review the work of \cite{Joyc14} on {\it deformation
theory} for compact SL $m$-folds with conical singularities.
Following \cite[Def.~5.4]{Joyc14}, we define the space $\M_\sX$
of compact SL $m$-folds $\hat X$ in $M$ with conical singularities
deforming a fixed SL $m$-fold $X$ with conical singularities.

\begin{dfn} Let $(M,J,\om,\Om)$ be an almost Calabi--Yau
$m$-fold and $X$ a compact SL $m$-fold in $M$ with conical
singularities at $x_1,\ldots,x_n$ with identifications
$\up_i:\C^m\ra T_{x_i}M$ and cones $C_1,\ldots,C_n$. Define
the {\it moduli space} $\M_\sX$ {\it of deformations of\/}
$X$ to be the set of $\hat X$ such that
\begin{itemize}
\setlength{\itemsep}{0pt}
\setlength{\parsep}{0pt}
\item[(i)] $\hat X$ is a compact SL $m$-fold in $M$ with
conical singularities at $\hat x_1,\ldots,\hat x_n$ with
cones $C_1,\ldots,C_n$, for some $\hat x_i$ and
identifications~$\hat\up_i:\C^m\ra T_{\smash{\hat x_i}}M$.
\item[(ii)] There exists a homeomorphism $\hat\io:X\ra\hat X$
with $\hat\io(x_i)=\hat x_i$ for $i=1,\ldots,n$ such that
$\hat\io\vert_{X'}:X'\ra\hat X'$ is a diffeomorphism and
$\hat\io$ and $\io$ are isotopic as continuous maps
$X\ra M$, where $\io:X\ra M$ is the inclusion.
\end{itemize}
\label{gr5def1}
\end{dfn}

In \cite[Def.~5.6]{Joyc14} we define a {\it topology} on
$\M_\sX$, and explain why it is the natural choice. We will
not repeat the complicated definition here; readers are
referred to \cite[\S 5]{Joyc14} for details. In
\cite[Th.~6.10]{Joyc14} we describe $\M_\sX$ near $X$, in terms
of a smooth map $\Phi$ between the {\it infinitesimal deformation
space} $\I_\sXp$ and the {\it obstruction space}~$\O_\sXp$.

\begin{thm} Suppose $(M,J,\om,\Om)$ is an almost Calabi--Yau
$m$-fold and\/ $X$ a compact SL\/ $m$-fold in $M$ with conical
singularities at\/ $x_1,\ldots,x_n$ and cones $C_1,\ldots,C_n$.
Let\/ $\M_\sX$ be the moduli space of deformations
of\/ $X$ as an SL\/ $m$-fold with conical singularities in $M$,
as in Definition \ref{gr5def1}. Set\/~$X'=X\sm\{x_1,\ldots,x_n\}$.

Then there exist natural finite-dimensional vector spaces
$\I_\sXp$, $\O_\sXp$ such that\/ $\I_\sXp$ is isomorphic to
the image of\/ $H^1_{\rm cs}(X',\R)$ in $H^1(X',\R)$ and\/ 
$\dim\O_\sXp=\sum_{i=1}^n\sind(C_i)$, where $\sind(C_i)$ is
the stability index of Definition \ref{gr3def3}. There exists
an open neighbourhood\/ $U$ of\/ $0$ in $\I_\sXp$, a smooth
map $\Phi:U\ra\O_\sXp$ with\/ $\Phi(0)=0$, and a map
$\Xi:\{u\in U:\Phi(u)=0\}\ra\M_\sX$ with\/ $\Xi(0)=X$ which is
a homeomorphism with an open neighbourhood of\/ $X$ in~$\M_\sX$.
\label{gr5thm1}
\end{thm}

Here is a sketch of the proof. For simplicity, consider
first the subset of $\hat X\in\M_\sX$ which have the
same singular points $x_1,\ldots,x_n$ and identifications
$\up_1,\ldots,\up_n$ as $X$. Generalizing Theorem
\ref{gr2thm3}, in \cite[Th.~4.3]{Joyc13} we define a
{\it Lagrangian neighbourhood\/} $U_\sXp,\Phi_\sXp$
for $X'$, with certain compatibilities with
$\Up_i,\phi_i$ near $x_i$. If $\hat X$ is $C^1$ close
to $X$ in an appropriate sense then $\hat X'=\Phi_\sXp
\bigl(\Ga(\al)\bigr)$, where $\Ga(\al)\subset U_\sXp$
is the graph of a small 1-form $\al$ on~$X'$.

Since $\hat X'$ is Lagrangian, $\al$ is {\it closed}, as
in \S\ref{gr21}. Also, applying Theorem \ref{gr4thm2} to
$X,\hat X$ and noting that $\al$ on $S_i$ corresponds to
$\hat\phi_i-\phi_i$ on $\Si_i\t(0,R')$, we find that if
$i=1,\ldots,n$ and $\mu_i'\in(2,3)$ with $(2,\mu_i']\cap
\D_\sSii=\emptyset$ then
\e
\bmd{\na^k\al(x)}=O\bigl(d(x,x_i)^{\mu_i'-1-k}\bigr)
\quad\text{near $x_i$ for all $k\ge 0$.}
\label{gr5eq1}
\e

As $\al$ is closed it has a cohomology class $[\al]\in H^1(X',\R)$.
Since \eq{gr5eq1} implies that $\al$ decays quickly near $x_i$, it
turns out that $\al$ must be {\it exact} near $x_i$. Therefore
$[\al]$ can be represented by a compactly-supported form on $X'$,
and lies in the image of $H^1_{\rm cs}(X',\R)$ in~$H^1(X',\R)$.

Choose a vector space $\I_\sXp$ of compactly-supported 1-forms
on $X'$ representing the image of $H^1_{\rm cs}(X',\R)$ in
$H^1(X',\R)$. Then we can write $\al=\be+\d f$, where $\be\in
\I_\sXp$ with $[\al]=[\be]$ is unique, and $f\in C^\iy(X')$ is
unique up to addition of constants. As $\hat X'$ is special
Lagrangian we find that $f$ satisfies a {\it second-order
nonlinear elliptic p.d.e.} similar to \eq{gr4eq3}:
\e
\d^*\bigl(\psi^m(\be+\d f)\bigr)(x)=
Q\bigl(x,(\be+\d f)(x),(\na\be+\na^2f)(x)\bigr)
\label{gr5eq2}
\e
for $x\in X'$. The {\it linearization} of \eq{gr5eq2} at
$\be=f=0$ is~$\d^*\bigl(\psi^m(\be+\d f)\bigr)=0$.

We study the family of small solutions $\be,f$ of \eq{gr5eq2}
for which $f$ has the decay near $x_i$ required by \eq{gr5eq1}.
There is a ready-made theory of analysis on manifolds with
cylindrical ends developed by Lockhart and McOwen \cite{Lock},
which is well-suited to this task. We work on certain {\it weighted
Sobolev spaces} $L^p_{k,{\bs\mu}}(X')$ of functions on~$X'$.

By results from \cite{Lock} it turns out that the
operator $f\mapsto\d^*(\psi^m\d f)$ is a {\it Fredholm} map
$L^p_{k,{\bs\mu}}(X')\ra L^p_{k-2,{\bs\mu}-2}(X')$, with
cokernel of dimension $\sum_{i=1}^nN_\sSii(2)$. This cokernel
is in effect the {\it obstruction space} to deforming $X$ with
$x_i,\up_i$ fixed, as it is the obstruction space to solving
the linearization of \eq{gr5eq2} in $f$ at~$\be=f=0$.

By varying the $x_i$ and $\up_i$, and allowing $f$ to converge
to different constant values on the ends of $X'$ rather than
zero, we can overcome many of these obstructions. This reduces
the dimension of the obstruction space $\O_\sXp$ from
$\sum_{i=1}^nN_\sSii(2)$ to $\sum_{i=1}^n\sind(C_i)$.  The
obstruction map $\Phi$ is constructed using the Implicit
Mapping Theorem for Banach spaces. This concludes our sketch.

If the $C_i$ are {\it stable} then $\O_\sXp=\{0\}$ and we
deduce~\cite[Cor.~6.11]{Joyc14}:

\begin{cor} Suppose $(M,J,\om,\Om)$ is an almost Calabi--Yau
$m$-fold and\/ $X$ a compact SL\/ $m$-fold in $M$ with stable
conical singularities, and let\/ $\M_\sX$ and\/ $\I_\sXp$ be
as in Theorem \ref{gr5thm1}. Then $\M_\sX$ is a smooth manifold
of dimension~$\dim\I_\sXp$.
\label{gr5cor1}
\end{cor}

This has clear similarities with Theorem \ref{gr2thm3}. Here
is another simple condition for $\M_\sX$ to be a manifold
near $X$, \cite[Def.~6.12]{Joyc14}.

\begin{dfn} Let $(M,J,\om,\Om)$ be an almost Calabi--Yau
$m$-fold and $X$ a compact SL $m$-fold in $M$ with conical
singularities, and let $\I_\sXp,\O_\sXp,U$ and $\Phi$ be as
in Theorem \ref{gr5thm1}. We call $X$ {\it transverse} if
the linear map $\d\Phi\vert_0:\I_\sXp\ra\O_\sXp$ is surjective.
\label{gr5def2}
\end{dfn}

If $X$ is transverse then $\{u\in U:\Phi(u)=0\}$ is a manifold
near 0, so Theorem \ref{gr5thm1} yields~\cite[Cor.~6.13]{Joyc14}:

\begin{cor} Suppose $(M,J,\om,\Om)$ is an almost Calabi--Yau $m$-fold
and\/ $X$ a transverse compact SL\/ $m$-fold in $M$ with conical
singularities, and let\/ $\M_\sX,\I_\sXp$ and\/ $\O_\sXp$ be as
in Theorem \ref{gr5thm1}. Then $\M_\sX$ is near $X$ a smooth manifold
of dimension~$\dim\I_\sXp-\dim\O_\sXp$.
\label{gr5cor2}
\end{cor}

Now there are a number of well-known moduli space problems in
geometry where in general moduli spaces are obstructed and
singular, but after a generic perturbation they become smooth
manifolds. For instance, moduli spaces of instantons on
4-manifolds can be made smooth by choosing a generic metric,
and similar things hold for Seiberg--Witten equations, and
moduli spaces of pseudo-holomorphic curves in symplectic
manifolds.

In \cite[\S 9]{Joyc14} we try (but do not quite succeed) to
replicate this for moduli spaces of SL $m$-folds with conical
singularities, by choosing a {\it generic K\"ahler metric} in
a fixed K\"ahler class. This is the idea
behind~\cite[Conj.~9.5]{Joyc14}:

\begin{conj} Let\/ $(M,J,\om,\Om)$ be an almost Calabi--Yau
$m$-fold, $X$ a compact SL\/ $m$-fold in $M$ with conical
singularities, and\/ $\I_\sXp,\O_\sXp$ be as in Theorem
\ref{gr5thm1}. Then for a second category subset of K\"ahler
forms $\check\om$ in the K\"ahler class of\/ $\om$, the moduli
space $\check\M_\sX$ of compact SL\/ $m$-folds $\hat X$ with
conical singularities in $(M,J,\check\om,\Om)$ isotopic to $X$
consists solely of transverse $\hat X$, and so is a manifold
of dimension~$\dim\I_\sXp-\dim\O_\sXp$.
\label{gr5conj}
\end{conj}

A partial proof of this is given in \cite[\S 9]{Joyc14}.
If we could treat the moduli spaces $\M_\sX$ as compact, the
conjecture would follow from \cite[Th.~9.3]{Joyc14}. However,
without knowing $\M_\sX$ is compact, the condition that
$\check\M_\sX$ is smooth everywhere is in effect the
intersection of an infinite number of genericity conditions
on $\check\om$, and we do not know that this intersection is
dense (or even nonempty) in the K\"ahler class.

Notice that Conjecture \ref{gr5conj} constrains the topology
and cones of SL $m$-folds $X$ with conical singularities that
can occur in a generic almost Calabi--Yau $m$-fold, as we
must have~$\dim\I_\sXp\ge\dim\O_\sXp$.

\section{Asymptotically Conical SL $m$-folds}
\label{gr6}

We now discuss {\it Asymptotically Conical\/} SL $m$-folds $L$
in $\C^m$, \cite[Def.~7.1]{Joyc13}.

\begin{dfn} Let $C$ be a closed SL cone in $\C^m$ with isolated
singularity at 0 for $m>2$, and let $\Si=C\cap{\cal S}^{2m-1}$,
so that $\Si$ is a compact, nonsingular $(m-1)$-manifold, not
necessarily connected. Let $g_\sSi$ be the metric on $\Si$
induced by the metric $g'$ on $\C^m$ in \eq{gr2eq1}, and $r$ the
radius function on $\C^m$. Define $\io:\Si\t(0,\iy)\ra\C^m$ by
$\io(\si,r)=r\si$. Then the image of $\io$ is $C\sm\{0\}$,
and $\io^*(g')=r^2g_\sSi+\d r^2$ is the cone metric
on~$C\sm\{0\}$.

Let $L$ be a closed, nonsingular SL $m$-fold in $\C^m$. We
call $L$ {\it Asymptotically Conical (AC)} with {\it rate}
$\la<2$ and {\it cone} $C$ if there exists a compact subset
$K\subset L$ and a diffeomorphism $\vp:\Si\t(T,\iy)\ra L\sm K$
for $T>0$, such that
\begin{equation*}
\bmd{\na^k(\vp-\io)}=O(r^{\la-1-k})
\quad\text{as $r\ra\iy$ for $k=0,1$.}
\end{equation*}
Here $\na,\md{\,.\,}$ are computed using the cone metric~$\io^*(g')$.
\label{gr6def1}
\end{dfn}

This is very similar to Definition \ref{gr3def4}, and in fact
there are strong parallels between the theories of SL $m$-folds
with conical singularities and of Asymptotically Conical SL
$m$-folds. We continue to assume $m>2$ throughout.

\subsection{Regularity and deformation theory of AC SL\/ $m$-folds}
\label{gr61}

Here are the analogues of Theorems \ref{gr4thm1} and \ref{gr4thm2},
proved in~\cite[Th.s~7.4 \& 7.11]{Joyc13}.

\begin{thm} Suppose $L$ is an AC SL\/ $m$-fold in $\C^m$ with cone
$C$, and let\/ $\Si,\io$ be as in Definition \ref{gr6def1}. Then
for sufficiently large $T>0$ there exist unique $K,\vp$ satisfying
the conditions of Definition \ref{gr6def1} and\/
$(\vp-\io)(\si,r)\perp T_{\io(\si,r)}C$ in $\C^m$ for
all\/~$(\si,r)\in\Si\t(T,\iy)$.
\label{gr6thm1}
\end{thm}

\begin{thm} In Theorem \ref{gr6thm1}, if either $\la=\la'$, or
$\la,\la'$ lie in the same connected component of\/ $\R\sm\D_\sSi$,
then $L$ is an AC SL\/ $m$-fold with rate $\la'$ and\/
$\bmd{\na^k(\vp-\io)}=O(r^{\la'-1-k})$ for all\/ $k\ge 0$. Here
$\na,\md{\,.\,}$ are computed using the cone metric $\io^*(g')$
on~$\Si\t(T,\iy)$.
\label{gr6thm2}
\end{thm}

The {\it deformation theory} of Asymptotically Conical SL
$m$-folds in $\C^m$ has been studied independently by Pacini
\cite{Paci} and Marshall \cite{Mars}. Pacini's results are
earlier, but Marshall's are more complete.

\begin{dfn} Suppose $L$ is an Asymptotically Conical SL
$m$-fold in $\C^m$ with cone $C$ and rate $\la<2$, as in
Definition \ref{gr6def1}. Define the {\it moduli space
$\M_\sL^\la$ of deformations of\/ $L$ with rate} $\la$
to be the set of AC SL $m$-folds $\hat L$ in $\C^m$ with
cone $C$ and rate $\la$, such that $\hat L$ is diffeomorphic
to $L$ and isotopic to $L$ as an Asymptotically Conical
submanifold of $\C^m$. One can define a natural {\it topology}
on $\M_\sL^\la$, in a similar way to the conical singularities
case of~\cite[Def.~5.6]{Joyc14}.
\label{gr6def2}
\end{dfn}

Note that if $L$ is an AC SL $m$-fold with rate $\la$, then
it is {\it also} an AC SL $m$-fold with rate $\la'$ for any
$\la'\in[\la,2)$. Thus we have defined a 1-{\it parameter
family} of moduli spaces $\M_\sL^{\smash{\la'}}$ for $L$, and
not just one. Since we did not impose any condition on $\la$
in Definition \ref{gr6def1} analogous to \eq{gr3eq8} in the
conical singularities case, it turns out that $\M_\sL^\la$
depends nontrivially on~$\la$.

The following result can be deduced from Marshall
\cite[Th.~6.2.15]{Mars} and \cite[Table~5.1]{Mars}.
(See also Pacini \cite[Th.~2 \& Th.~3]{Paci}.) It implies
conjectures by the author in \cite[Conj.~2.12]{Joyc1}
and~\cite[\S 10.2]{Joyc8}.

\begin{thm} Let\/ $L$ be an Asymptotically Conical SL\/
$m$-fold in $\C^m$ with cone $C$ and rate $\la<2$, and
let\/ $\M_\sL^\la$ be as in Definition \ref{gr6def2}.
Set\/ $\Si=C\cap{\cal S}^{2m-1}$, and let\/ $\D_\sSi,
N_\sSi$ be as in \S\ref{gr31} and\/ $b^k(L),b^k_{\rm
cs}(L)$ as in \S\ref{gr34}. Then
\begin{itemize}
\setlength{\itemsep}{0pt}
\setlength{\parsep}{0pt}
\item[{\rm(a)}] If\/ $\la\in(0,2)\sm\D_\sSi$ then
$\M_\sL^\la$ is a manifold with
\begin{equation*}
\dim\M_\sL^\la=b^1(L)-b^0(L)+N_\sSi(\la).
\end{equation*}
Note that if\/ $0<\la<\min\bigl(\D_\sSi\cap
(0,\iy)\bigr)$ then~$N_\sSi(\la)=b^0(\Si)$.
\item[{\rm(b)}] If\/ $\la\in(2-m,0)$ then $\M_\sL^\la$
is a manifold of dimension~$b^1_{\rm cs}(L)=b^{m-1}(L)$.
\end{itemize}
\label{gr6thm3}
\end{thm}

This is the analogue of Theorems \ref{gr2thm3} and \ref{gr5thm1}
for AC SL $m$-folds. If $\la\in(2-m,2)\sm\D_\sSi$ then the
deformation theory for $L$ with rate $\la$ is {\it unobstructed\/}
and $\M_\sL^\la$ is a {\it smooth manifold\/} with a given
dimension. This is similar to the case of nonsingular compact
SL $m$-folds in Theorem \ref{gr2thm3}, but different to the
conical singularities case in Theorem~\ref{gr5thm1}.

\subsection{Cohomological invariants of AC SL $m$-folds}
\label{gr62}

Let $L$ be an AC SL $m$-fold in $\C^m$ with cone $C$, and set
$\Si=C\cap{\cal S}^{2m-1}$. Using the notation of \S\ref{gr34},
as in \eq{gr3eq10} there is a long exact sequence
\e
\cdots\ra
H^k_{\rm cs}(L,\R)\ra H^k(L,\R)\ra H^k(\Si,\R)\ra
H^{k+1}_{\rm cs}(L,\R)\ra\cdots.
\label{gr6eq1}
\e
Following \cite[Def.~7.2]{Joyc13} we define {\it cohomological
invariants\/} $Y(L),Z(L)$ of~$L$.

\begin{dfn} Let $L$ be an AC SL $m$-fold in $\C^m$ with cone $C$,
and let $\Si=C\cap{\cal S}^{2m-1}$. As $\om',\Im\Om'$ in \eq{gr2eq1}
are closed forms with $\om'\vert_L\equiv\Im\Om'\vert_L\equiv 0$, they
define classes in the relative de Rham cohomology groups $H^k(\C^m;
L,\R)$ for $k=2,m$. But for $k>1$ we have the exact sequence
\begin{equation*}
0=H^{k-1}(\C^m,\R)\ra H^{k-1}(L,\R){\buildrel\cong\over\longra}
H^k(\C^m;L,\R)\ra H^k(\C^m,\R)=0.
\end{equation*}
Let $Y(L)\in H^1(\Si,\R)$ be the image of $[\om']$ in
$H^2(\C^m;L,\R)\cong H^1(L,\R)$ under $H^1(L,\R)\ra H^1(\Si,R)$
in \eq{gr6eq1}, and $Z(L)\in H^{m-1}(\Si,\R)$ be the image of
$[\Im\Om']$ in $H^m(\C^m;L,\R)\cong H^{m-1}(L,\R)$ under
$H^{m-1}(L,\R)\ra H^{m-1}(\Si,R)$ in~\eq{gr6eq1}.
\label{gr6def3}
\end{dfn}

Here are some conditions for $Y(L)$ or $Z(L)$ to be
zero,~\cite[Prop.~7.3]{Joyc13}.

\begin{prop} Let\/ $L$ be an AC SL\/ $m$-fold in $\C^m$ with
cone $C$ and rate $\la$, and let\/ $\Si=C\cap{\cal S}^{2m-1}$.
If\/ $\la<0$ or $b^1(L)=0$ then $Y(L)=0$. If\/ $\la<2-m$ or
$b^0(\Si)=1$ then~$Z(L)=0$.
\label{gr6prop}
\end{prop}

\subsection{Examples}
\label{gr63}

Examples of AC SL $m$-folds $L$ are constructed by Harvey and Lawson
\cite[\S III.3]{HaLa}, the author \cite{Joyc2,Joyc3,Joyc4,Joyc6},
and others. Nearly all the known examples (up to translations) have
minimum rate $\la$ either 0 or $2-m$, which are topologically
significant values by Proposition \ref{gr6prop}. For instance, all
examples in \cite{Joyc3} have $\la=0$, and \cite[Th.~6.4]{Joyc2}
constructs AC SL $m$-folds with $\la=2-m$ in $\C^m$ from any SL
cone $C$ in $\C^m$. The only explicit, nontrivial examples known
to the author with $\la\ne 0,2-m$ are in \cite[Th.~11.6]{Joyc4},
and have~$\la=\frac{3}{2}$.

We shall give three families of examples of AC SL $m$-folds $L$
in $\C^m$ explicitly. The first family is adapted from Harvey
and Lawson~\cite[\S III.3.A]{HaLa}.

\begin{ex} Let $C_{\sst\rm HL}^m$ be the SL cone in $\C^m$ of
Example \ref{gr3ex1}. We shall define a family of AC SL $m$-folds
in $\C^m$ with cone $C_{\sst\rm HL}^m$. Let $a_1,\ldots,a_m\ge 0$
with exactly two of the $a_j$ zero and the rest positive. Write
${\bf a}=(a_1,\ldots,a_m)$, and define
\e
\begin{split}
L_{\sst\rm HL}^{\bf a}=\bigl\{&(z_1,\ldots,z_m)\in\C^m:
i^{m+1}z_1\cdots z_m\in[0,\iy),\\
&\ms{z_1}-a_1=\cdots=\ms{z_m}-a_m\bigr\}.
\end{split}
\label{gr6eq2}
\e
Then $L_{\sst\rm HL}^{\bf a}$ is an AC SL $m$-fold in $\C^m$
diffeomorphic to $T^{m-2}\t\R^2$, with cone $C_{\sst\rm HL}^m$
and rate 0. It is invariant under the $\U(1)^{m-1}$ group
\eq{gr3eq5}. It is surprising that equations of the form
\eq{gr6eq2} should define a nonsingular submanifold of
$\C^m$ {\it without boundary}, but in fact they do.

Now suppose for simplicity that $a_1,\ldots,a_{m-2}>0$ and
$a_{m-1}=a_m=0$. As $\Si_{\sst\rm HL}^m\cong T^{m-1}$ we have
$H^1(\Si_{\sst\rm HL}^m,\R)\cong\R^{m-1}$, and calculation
shows that $Y(L_{\sst\rm HL}^{\bf a})=(\pi a_1,\ldots,\pi
a_{m-2},0)\in\R^{m-1}$ in the natural coordinates. Since
$L_{\sst\rm HL}^{\bf a}\cong T^{m-2}\t\R^2$ we have
$H^1(L_{\sst\rm HL}^{\bf a},\R)=\R^{m-2}$, and
$Y(L_{\sst\rm HL}^{\bf a})$ lies in the image $\R^{m-2}
\subset\R^{m-1}$ of $H^1(L_{\sst\rm HL}^{\bf a},\R)$ in
$H^1(\Si_{\sst\rm HL}^m,\R)$, as in Definition
\ref{gr6def3}. As $b^0(\Si_{\sst\rm HL}^m)=1$, Proposition
\ref{gr6prop} shows that~$Z(L_{\sst\rm HL}^{\bf a})=0$.

Take $C=C_{\sst\rm HL}^m$, $\Si=\Si_{\sst\rm HL}^m$ and
$L=L_{\sst\rm HL}^{\bf a}$ in Theorem \ref{gr6thm3},
and let $0<\la<\min\bigl(\D_\sSi\cap(0,\iy)\bigr)$. Then
$b^1(L)=m-2$, $b^0(L)=1$ and $N_\sSi(\la)=b^0(\Si)=1$, so part
(a) of Theorem \ref{gr6thm3} shows that $\dim\M_\sL^\la=m-2$.
This is consistent with the fact that $L$ depends on $m-2$
real parameters~$a_1,\ldots,a_{m-2}>0$.

The family of all $L_{\sst\rm HL}^{\bf a}$ has $\ha m(m-1)$
connected components, indexed by which two of $a_1,\ldots,a_m$
are zero. Using the theory of \S\ref{gr7}, these can give many
{\it topologically distinct\/} ways to desingularize SL
$m$-folds with conical singularities with these cones.
\label{gr6ex1}
\end{ex}

Our second family, from \cite[Ex.~9.4]{Joyc2}, was chosen as
it is easy to write down.

\begin{ex} Let $m,a_1,\ldots,a_m,k$ and $L^{a_1,\ldots,a_m}_0$
be as in Example \ref{gr3ex2}. For $0\ne c\in\R$ define
\begin{align*}
L^{a_1,\ldots,a_m}_c=\bigl\{
\bigl(i{\rm e}^{ia_1\th}x_1&,{\rm e}^{ia_2\th}x_2,\ldots,
{\rm e}^{ia_m\th}x_m\bigr):
\th\in[0,2\pi),\\ 
&x_1,\ldots,x_m\in\R,\qquad 
a_1x_1^2+\cdots+a_mx_m^2=c\bigr\}.
\end{align*}
Then $L^{a_1,\ldots,a_m}_c$ is an AC SL $m$-fold in $\C^m$ with
rate 0 and cone $L^{a_1,\ldots,a_m}_0$. It is diffeomorphic as an
immersed SL $m$-fold to $({\cal S}^{k-1}\t\R^{m-k}\t{\cal S}^1)/\Z_2$
if $c>0$, and to $(\R^k\t{\cal S}^{m-k-1}\t{\cal S}^1)/\Z_2$ if~$c<0$.
\label{gr6ex2}
\end{ex}

Our third family was first found by Lawlor \cite{Lawl}, made more
explicit by Harvey \cite[p.~139--140]{Harv}, and discussed from a
different point of view by the author in \cite[\S 5.4(b)]{Joyc3}.
Our treatment is based on that of Harvey.

\begin{ex} Let $m>2$ and $a_1,\ldots,a_m>0$, and define
polynomials $p,P$ by
\begin{equation*}
p(x)=(1+a_1x^2)\cdots(1+a_mx^2)-1
\quad\text{and}\quad P(x)=\frac{p(x)}{x^2}.
\end{equation*}
Define real numbers $\phi_1,\ldots,\phi_m$ and $A$ by
\e
\phi_k=a_k\int_{-\iy}^\iy\frac{\d x}{(1+a_kx^2)\sqrt{P(x)}}
\quad\text{and}\quad A=\om_m(a_1\cdots a_m)^{-1/2},
\label{gr6eq3}
\e
where $\om_m$ is the volume of the unit sphere in $\R^m$. Clearly
$\phi_k,A>0$. But writing $\phi_1+\cdots+\phi_m$ as one integral gives
\begin{equation*}
\phi_1+\cdots+\phi_m=\int_0^\iy\frac{p'(x)\d x}{(p(x)+1)\sqrt{p(x)}}
=2\int_0^\iy\frac{\d w}{w^2+1}=\pi,
\end{equation*}
making the substitution $w=\sqrt{p(x)}$. So $\phi_k\in(0,\pi)$
and $\phi_1+\cdots+\phi_m=\pi$. This yields a 1-1 correspondence
between $m$-tuples $(a_1,\ldots,a_m)$ with $a_k>0$, and
$(m\!+\!1)$-tuples $(\phi_1,\ldots,\phi_m,A)$ with $\phi_k\in
(0,\pi)$, $\phi_1+\cdots+\phi_m=\pi$ and~$A>0$.

For $k=1,\ldots,m$ and $y\in\R$, define a function $z_k:\R\ra\C$ by
\begin{equation*}
z_k(y)={\rm e}^{i\psi_k(y)}\sqrt{a_k^{-1}+y^2}, \quad\text{where}\quad
\psi_k(y)=a_k\int_{-\iy}^y\frac{\d x}{(1+a_kx^2)\sqrt{P(x)}}\,.
\end{equation*}
Now write ${\bs\phi}=(\phi_1,\ldots,\phi_n)$, and define 
a submanifold $L^{{\bs\phi},A}$ in $\C^m$ by
\begin{equation*}
L^{{\bs\phi},A}=\bigl\{(z_1(y)x_1,\ldots,z_m(y)x_m):
y\in\R,\; x_k\in\R,\; x_1^2+\cdots+x_m^2=1\bigr\}.
\end{equation*}

Then $L^{{\bs\phi},A}$ is closed, embedded, and diffeomorphic
to ${\cal S}^{m-1}\t\R$, and Harvey \cite[Th.~7.78]{Harv} shows
that $L^{{\bs\phi},A}$ is {\it special Lagrangian}. One can also
show that $L^{{\bs\phi},A}$ is {\it Asymptotically Conical},
with rate $2-m$ and cone the union $\Pi^0\cup\Pi^{\bs\phi}$ of
two special Lagrangian $m$-planes $\Pi^0,\Pi^{\bs\phi}$ in
$\C^m$ given by
\begin{equation*}
\Pi^0=\bigl\{(x_1,\ldots,x_m):x_j\in\R\bigr\}
\;\>\text{and}\;\>
\Pi^{\bs\phi}=\bigl\{({\rm e}^{i\phi_1}x_1,\ldots,
{\rm e}^{i\phi_m}x_m):x_j\in\R\bigr\}.
\end{equation*}

As $\la=2-m<0$ we have $Y(L^{{\bs\phi},A})=0$ by Proposition
\ref{gr6prop}. Now $L^{{\bs\phi},A}\cong{\cal S}^{m-1}\t\R$
so that $H^{m-1}(L^{{\bs\phi},A},\R)\cong\R$, and $\Si=(\Pi^0
\cup\Pi^{\bs\phi})\cap{\cal S}^{2m-1}$ is the disjoint union
of two unit $(m\!-\!1)$-spheres ${\cal S}^{m-1}$, so $H^{m-1}
(\Si,\R)\cong\R^2$. The image of $H^{m-1}(L^{{\bs\phi},A},\R)$
in $H^{m-1}(\Si,\R)$ is $\bigl\{(x,-x):x\in\R\bigr\}$ in the
natural coordinates. Calculation shows that $Z(L^{{\bs\phi},A})=
(A,-A)\in H^{m-1}(\Si,\R)$, which is why we defined $A$ this way
in~\eq{gr6eq3}.

Apply Theorem \ref{gr6thm3} with $L=L^{{\bs\phi},A}$ and
$\la\in(2-m,0)$. As $L\cong{\cal S}^{m-1}\t\R$ we have
$b^1_{\rm cs}(L)=1$, so part (b) of Theorem \ref{gr6thm3} shows
that $\dim\M_\sL^\la=1$. This is consistent with the fact that
when $\bs\phi$ is fixed, $L^{{\bs\phi},A}$ depends on one real
parameter $A>0$. Here $\bs\phi$ is fixed in $\M_\sL^\la$ as the
cone $C=\Pi^0\cup\Pi^{\bs\phi}$ of $L$ depends on $\bs\phi$, and
all $\hat L\in\M_\sL^\la$ have the same cone $C$, by definition.
\label{gr6ex3}
\end{ex}

\section{Desingularizing singular SL $m$-folds}
\label{gr7}

We now discuss the work of \cite{Joyc15,Joyc16} on
{\it desingularizing} compact SL $m$-folds with conical
singularities. Here is the basic idea. Let $(M,J,\om,\Om)$
be an almost Calabi--Yau $m$-fold, and $X$ a compact SL
$m$-fold in $M$ with conical singularities $x_1,\ldots,x_n$
and cones $C_1,\ldots,C_n$. Suppose $L_1,\ldots,L_n$ are AC SL
$m$-folds in $\C^m$ with the same cones $C_1,\ldots,C_n$ as~$X$.

If $t>0$ then $tL_i=\{t\,{\bf x}:{\bf x}\in L_i\}$ is also
an AC SL $m$-fold with cone $C_i$. We construct a 1-parameter
family of compact, nonsingular {\it Lagrangian} $m$-folds $N^t$
in $(M,\om)$ for $t\in(0,\de)$ by gluing $tL_i$ into $X$ at
$x_i$, using a partition of unity.

When $t$ is small, $N^t$ is {\it close to special Lagrangian}
(its phase is nearly constant), but also {\it close to singular}
(it has large curvature and small injectivity radius). We prove
using analysis that for small $t\in(0,\de)$ we can deform $N^t$
to a {\it special\/} Lagrangian $m$-fold $\smash{\ti N^t}$ in
$M$, using a small Hamiltonian deformation.

The proof involves a delicate balancing act, showing that the
advantage of being close to special Lagrangian outweighs the
disadvantage ofbeing nearly singular. Doing this in full
generality is rather complex. Here is our simplest
desingularization result,~\cite[Th.~6.13]{Joyc15}.

\begin{thm} Suppose $(M,J,\om,\Om)$ is an almost Calabi--Yau
$m$-fold and\/ $X$ a compact SL\/ $m$-fold in $M$ with conical
singularities at\/ $x_1,\ldots,x_n$ and cones $C_1,\ldots,C_n$.
Let\/ $L_1,\ldots,L_n$ be Asymptotically Conical SL\/ $m$-folds
in $\C^m$ with cones $C_1,\ldots,C_n$ and rates $\la_1,\ldots,
\la_n$. Suppose $\la_i<0$ for $i=1,\ldots,n$, and\/ $X'=X\sm
\{x_1,\ldots,x_n\}$ is connected.

Then there exists $\ep>0$ and a smooth family $\bigl\{
\smash{\ti N^t}:t\in(0,\ep]\bigr\}$ of compact, nonsingular
SL\/ $m$-folds in $(M,J,\om,\Om)$, such that\/ $\smash{\ti N^t}$
is constructed by gluing $tL_i$ into $X$ at\/ $x_i$ for
$i=1,\ldots,n$. In the sense of currents, $\smash{\ti N^t}\ra
X$ as~$t\ra 0$.
\label{gr7thm1}
\end{thm}

The theorem contains two {\it simplifying assumptions}:
\begin{itemize}
\setlength{\itemsep}{0pt}
\setlength{\parsep}{0pt}
\item[(a)] that $X'$ is connected, and
\item[(b)] that $\la_i<0$ for all~$i$.
\end{itemize}
These avoid two kinds of {\it obstructions} to desingularizing
$X$ using the~$L_i$.

In \cite[Th.~7.10]{Joyc15} we remove assumption (a), allowing
$X'$ not connected.

\begin{thm} Suppose $(M,J,\om,\Om)$ is an almost Calabi--Yau
$m$-fold and\/ $X$ a compact SL\/ $m$-fold in $M$ with conical
singularities at\/ $x_1,\ldots,x_n$ and cones $C_1,\ldots,C_n$.
Define $\psi:M\ra(0,\iy)$ as in \eq{gr2eq3}. Let\/ $L_1,\ldots,L_n$
be Asymptotically Conical SL\/ $m$-folds in $\C^m$ with cones
$C_1,\ldots,C_n$ and rates $\la_1,\ldots,\la_n$. Suppose $\la_i<0$
for $i=1,\ldots,n$. Write $X'=X\sm\{x_1,\ldots,x_n\}$
and\/~$\Si_i=C_i\cap{\cal S}^{2m-1}$.

Set\/ $q=b^0(X')$, and let\/ $X_1',\ldots,X_q'$ be the connected
components of\/ $X'$. For $i=1,\ldots,n$ let\/ $l_i=b^0(\Si_i)$,
and let\/ $\Si_i^1,\ldots,\Si_i^{\smash{l_i}}$ be the connected
components of\/ $\Si_i$. Define $k(i,j)=1,\ldots,q$ by $\Up_i
\circ\vp_i\bigl(\Si_i^j\t(0,R')\bigr)\subset X'_{\smash{k(i,j)}}$
for $i=1,\ldots,n$ and $j=1,\ldots,l_i$. Suppose that
\e
\sum_{\substack{1\le i\le n, \; 1\le j\le l_i: \\
k(i,j)=k}}\psi(x_i)^mZ(L_i)\cdot[\Si_i^j\,]=0
\quad\text{for all\/ $k=1,\ldots,q$.}
\label{gr7eq1}
\e

Suppose also that the compact\/ $m$-manifold\/ $N$ obtained by
gluing $L_i$ into $X'$ at\/ $x_i$ for $i=1,\ldots,n$ is connected.
A sufficient condition for this to hold is that\/ $X$ and\/ $L_i$
for $i=1,\ldots,n$ are connected.

Then there exists $\ep>0$ and a smooth family $\smash{\bigl\{
\ti N^t:t\in(0,\ep]\bigr\}}$ of compact, nonsingular SL\/
$m$-folds in $(M,J,\om,\Om)$ diffeomorphic to $N$, such
that\/ $\smash{\ti N^t}$ is constructed by gluing $tL_i$
into $X$ at\/ $x_i$ for $i=1,\ldots,n$. In the sense of
currents in Geometric Measure Theory, $\smash{\ti N^t}\ra
X$ as~$t\ra 0$.
\label{gr7thm2}
\end{thm}

The new issue here is that if $X'$ is not connected then there
is an {\it analytic obstruction} to deforming $N^t$ to
$\smash{\ti N^t}$, because the Laplacian $\De^t$ on functions
on $N^t$ has {\it small eigenvalues} of size $O(t^{m-2})$. As
in \S\ref{gr62} the $L_i$ have {\it cohomological invariants}
$Z(L_i)$ in $H^{m-1}(\Si_i,\R)$ derived from the relative
cohomology class of $\Im\Om'$. It turns out that we can only
deform $N^t$ to $\smash{\ti N^t}$ if the $Z(L_i)$ satisfy
\eq{gr7eq1}. This equation arises by requiring the projection
of an error term to the eigenspaces of $\De^t$ with small
eigenvalues to be zero.

In \cite[Th.~6.13]{Joyc16} we remove assumption (b), extending Theorem
\ref{gr7thm1} to the case $\la_i\le 0$, and allowing~$Y(L_i)\ne 0$.

\begin{thm} Let\/ $(M,J,\om,\Om)$ be an almost Calabi--Yau $m$-fold
for $2\!<\!m\!<\nobreak\!6$, and\/ $X$ a compact SL\/ $m$-fold
in $M$ with conical singularities at\/ $x_1,\ldots,x_n$ and cones
$C_1,\ldots,C_n$. Let\/ $L_1,\ldots,L_n$ be Asymptotically
Conical SL\/ $m$-folds in $\C^m$ with cones $C_1,\ldots,C_n$
and rates $\la_1,\ldots,\la_n$. Suppose that\/ $\la_i\le 0$
for\/ $i=1,\ldots,n$, that\/ $X'=X\sm\{x_1,\ldots,x_n\}$ is
connected, and that there exists $\varrho\in H^1(X',\R)$
such that\/ $\bigl(Y(L_1),\ldots,Y(L_n)\bigr)$ is the image
of\/ $\varrho$ under the map $H^1(X',\R)\ra\bigoplus_{i=1}^n
H^1(\Si_i,\R)$ in \eq{gr3eq10}, where~$\Si_i=C_i\cap{\cal S}^{2m-1}$.

Then there exists $\ep>0$ and a smooth family $\bigl\{
\smash{\ti N^t}:t\in(0,\ep]\bigr\}$ of compact, nonsingular SL\/
$m$-folds in $(M,J,\om,\Om)$, such that\/ $\smash{\ti N^t}$ is
constructed by gluing $tL_i$ into $X$ at\/ $x_i$ for
$i=1,\ldots,n$. In the sense of currents, $\smash{\ti N^t}\ra
X$ as~$t\ra 0$.
\label{gr7thm3}
\end{thm}

From \S\ref{gr63}, the $L_i$ have {\it cohomological invariants}
$Y(L_i)$ in $H^1(\Si_i,\R)$ derived from the relative cohomology
class of $\om'$. The new issue in Theorem \ref{gr7thm3} is that
if $Y(L_i)\ne 0$ then there are obstructions to the existence of
$N^t$ as a {\it Lagrangian} $m$-fold. That is, we can only define
$N^t$ if the $Y(L_i)$ satisfy an equation. This did not appear in
Theorem \ref{gr7thm1}, as $\la_i<0$ implies that~$Y(L_i)=0$.

To define the $N^t$ when $Y(L_i)\ne 0$ we must also use a
more complicated construction. This introduces new errors. To
overcome these errors when we deform $N^t$ to $\smash{\ti N^t}$
we must assume that $m<6$. There is also \cite[Th.~6.12]{Joyc16}
a result combining the modifications of Theorems \ref{gr7thm2}
and \ref{gr7thm3}, but for brevity we will not give it.

\section{Directions for future research}
\label{gr8}

Finally we discuss directions the field of special Lagrangian
singularities might develop in the future, giving a number
of problems the author believes are worth attention. Some
of these problems may be too difficult to solve completely,
but can still serve as a guide.

\subsection{The index of singularities of SL\/ $m$-folds}
\label{gr81}

We now consider the {\it boundary\/} $\pd\M_\sN$ of a
moduli space $\M_\sN$ of SL $m$-folds.

\begin{dfn} Let $(M,J,\om,\Om)$ be an almost Calabi--Yau
$m$-fold, $N$ a compact, nonsingular SL $m$-fold in $M$,
and $\M_\sN$ the moduli space of deformations of $N$ in $M$.
Then $\M_\sN$ is a smooth manifold of dimension $b^1(N)$,
in general noncompact. We can construct a natural {\it
compactification} $\oM_\sN$ as follows.

Regard $\M_\sN$ as a moduli space of special Lagrangian
{\it integral currents} in the sense of Geometric Measure
Theory, as discussed in \cite[\S 6]{Joyc13}. Let $\oM_\sN$
be the closure of $\M_\sN$ in the space of integral currents.
As elements of $\M_\sN$ have uniformly bounded volume,
$\oM_\sN$ is {\it compact}. Define the {\it boundary}
$\pd\M_\sN$ to be $\oM_\sN\sm\M_\sN$. Then elements of
$\pd\M_\sN$ are {\it singular SL integral currents}.
\label{gr8def1}
\end{dfn}

In good cases, say if $(M,J,\om,\Om)$ is suitably generic,
it seems reasonable that $\pd\M_\sN$ should be divided into
a number of {\it strata}, each of which is a moduli space of
singular SL $m$-folds with singularities of a particular type,
and is itself a manifold with singularities. In particular, some
or all of these strata could be moduli spaces $\M_\sX$ of SL
$m$-folds with isolated conical singularities, as in~\S\ref{gr5}. 

Suppose $\M_\sN$ is a moduli space of compact, nonsingular SL
$m$-folds $N$ in $(M,J,\om,\Om)$, and $\M_\sX$ a moduli space
of singular SL $m$-folds in $\pd\M_\sN$ with singularities
of a particular type, and $X\in\M_\sX$. Following
\cite[\S 8.3]{Joyc17}, we (loosely) define the {\it index} of
the singularities of $X$ to be $\ind(X)=\dim\M_\sN-\dim\M_\sX$,
provided $\M_\sX$ is smooth near $X$. Note that $\ind(X)$
depends on $N$ as well as~$X$.

In \cite[Th.~8.10]{Joyc17} we use the results of
\cite{Joyc14,Joyc15,Joyc16} to compute $\ind(X)$ when $X$
is {\it transverse} with conical singularities, in the sense
of Definition \ref{gr5def2}. Here is a simplified version of
the result, where we assume that $H^1_{\rm cs}(L_i,\R)\ra
H^1(L_i,\R)$ is surjective to avoid a complicated correction
term to $\ind(X)$ related to the obstructions to defining
$N^t$ as a Lagrangian $m$-fold.

\begin{thm} Let\/ $X$ be a compact, transverse SL $m$-fold in
$(M,J,\om,\Om)$ with conical singularities at\/ $x_1,\ldots,x_n$
and cones $C_1,\ldots,C_n$. Let\/ $L_1,\ldots,L_n$ be AC SL\/
$m$-folds in $\C^m$ with cones $C_1,\ldots,C_n$, such that the
natural projection $H^1_{\rm cs}(L_i,\R)\ra H^1(L_i,\R)$ is
surjective. Construct desingularizations $N$ of\/ $X$ by gluing
AC SL\/ $m$-folds $L_1,\ldots,L_n$ in at\/ $x_1,\ldots,x_n$, as
in \S\ref{gr7}. Then
\e
\ind(X)=1-b^0(X')+\ts\sum_{i=1}^nb^1_{\rm cs}(L_i)+
\sum_{i=1}^n\sind(C_i).
\label{gr8eq1}
\e
\label{gr8thm1}
\end{thm}

If the cones $C_i$ are not {\it rigid\/}, for instance if
$C_i\sm\{0\}$ is not connected, then \eq{gr8eq1} should
be corrected, as in \cite[\S 8.3]{Joyc17}. If Conjecture
\ref{gr5conj} is true then for a {\it generic} K\"ahler
form $\om$, {\it all\/} compact SL $m$-folds $X$ with
conical singularities are transverse, and so Theorem
\ref{gr8thm1} and \cite[Th.~8.10]{Joyc17} allow us to
calculate~$\ind(X)$. 

Now singularities with {\it small index} are the most commonly
occurring, and so arguably the most interesting kinds of
singularity. Also, as $\ind(X)\le\dim\M_\sN$, for various
problems, such as those in \S\ref{gr83} and \S\ref{gr84}, it
will only be necessary to know about singularities with index
up to a certain value. This motivates the following:

\begin{prob} Classify types of singularities of SL 3-folds
with {\it small index} in suitably generic almost Calabi--Yau
3-folds, say with index 1,2 or~3.
\label{gr8prob1}
\end{prob}

Here we restrict to $m=3$ to make the problem more feasible,
though still difficult. Note, however, that we do {\it not\/}
restrict to isolated conical singularities, so a complete,
rigorous answer would require a theory of more general kinds
of singularities of SL 3-folds.

One can make some progress on this problem simply by studying
the many examples of singular SL 3-folds in \cite{HaLa,Hask1}
and \cite{Joyc3,Joyc4,Joyc5,Joyc6,Joyc7,Joyc8,Joyc9,Joyc10,
Joyc11,Joyc12}, calculating or guessing the index of each, and
ruling out other kinds of singularities by plausible-sounding
arguments. Using these techniques I have a conjectural
classification of index 1 singularities of SL 3-folds, which
involves the SL $T^2$-cone $C_{\sst\rm HL}^3$ of \eq{gr3eq4},
and several different kinds of singularity whose tangent cone
is two copies of $\R^3$, intersecting in 0, $\R$ or~$\R^3$.

Coming from another direction, {\it integrable systems}
techniques may yield rigorous classification results for SL
$T^2$-cones by index. Haskins \cite[Th.~A]{Hask2} has used
them to prove that the SL $T^2$-cone $C_{\sst\rm HL}^3$ in
$\C^3$ of \eq{gr3eq4} is up to $\SU(3)$ equivalence the
{\it unique} SL $T^2$-cone $C$ with $\sind(C)=0$. Now the
index of a singularity modelled on $C$ is at least
$\sind(C)+1$, so this implies that $C_{\sst\rm HL}^3$ is the
unique SL $T^2$-cone with index 1 in Problem~\ref{gr8prob1}.

\subsection{Singularities which are not isolated conical}
\label{gr82}

Singularities of SL $m$-folds which are not `isolated conical
singularities' in the sense of Definition \ref{gr3def4} are
an important, but virtually unexplored, subject. Here are
some known classes of nontrivial examples when~$m=3$.
\begin{itemize}
\setlength{\itemsep}{0pt}
\setlength{\parsep}{0pt}
\item[(i)] In \cite{Joyc6} we study {\it ruled\/} SL 3-folds in
$\C^3$, that is, SL 3-folds $N$ fibred by a 2-dimensional
family $\Si$ of real straight lines in $\C^3$. When $\Si$
is nonsingular $N$ can still have singularities, and examples
may be written down very explicitly, as in~\cite[Th.~7.1]{Joyc6}.

The {\it tangent cones} of such singularities, in the
sense of Geometric Measure Theory, are generally $\R^3$
with multiplicity $k>1$. Near the singular point, the SL
3-fold resembles a $k$-fold branched cover of $\R^3$,
branched along $\R$. A similar class of singularities
of SL 3-folds, with tangent cone $\R^3$ with multiplicity
2, is studied in~\cite[\S 6]{Joyc4}.
\item[(ii)] In \cite{Joyc9,Joyc10,Joyc11} we study SL 3-folds
in $\C^3$ invariant under the $\U(1)$-action
\begin{equation*}
{\rm e}^{i\th}:(z_1,z_2,z_3)\mapsto
({\rm e}^{i\th}z_1,{\rm e}^{-i\th}z_2,z_3)
\quad\text{for ${\rm e}^{i\th}\in\U(1)$.}
\end{equation*}
The three papers are surveyed in \cite{Joyc12}. A $\U(1)$-invariant
SL 3-fold $N$ may locally be written in the form
\begin{align*}
N=\bigl\{(z_1,z_2,z_3)\in\C^3:\,& z_1z_2=v(x,y)+iy,\quad z_3=x+iu(x,y),\\
&\ms{z_1}-\ms{z_2}=2a,\quad (x,y)\in S\bigr\},
\end{align*}
where $S$ is a domain in $\R^2$, $a\in\R$ and $u,v:S\ra\R$
satisfy (in a weak sense if $a=0$) the {\it nonlinear
Cauchy--Riemann equations}
\e
\frac{\pd u}{\pd x}=\frac{\pd v}{\pd y}\quad\text{and}\quad
\frac{\pd v}{\pd x}=-2\bigl(v^2+y^2+a^2\bigr)^{1/2}\frac{\pd u}{\pd y}.
\label{gr8eq2}
\e

Using analytic techniques, we construct and study solutions
$u,v$ of \eq{gr8eq2} satisfying boundary conditions on a
strictly convex domain $S$. These include many {\it singular
solutions}, and we show in \cite[\S 9--\S 10]{Joyc11} that
we can construct  countably many distinct geometrical-topological
types of isolated SL 3-fold singularities, whose tangent cone is
the union of two $\R^3$'s in $\C^3$, intersecting in~$\R$.
\end{itemize}

There appear to the author to be two ways of studying
special Lagrangian singularities which are not isolated
conical. The first is to try and study {\it all\/}
singularities of special Lagrangian integral currents,
using Geometric Measure Theory. As far as the author
understands (which is not very far), it will be difficult
to use the special Lagrangian condition in GMT, or to say
anything nontrivial about special Lagrangian singularities
in this generality.

The second way is to define some restricted class of singularities
and then study them, just as we did in \S\ref{gr3}--\S\ref{gr7}.
The problem here is to decide upon a suitable kind of {\it local
model\/} for the singularities, and appropriate {\it asymptotic
conditions} for how the SL $m$-fold approaches the local model
near the singularity. Now not just any local model and asymptotic
conditions will do.

For a class of singularities to be worth studying, they should
occur reasonably often in `real life', so that, for instance,
examples of such singularities might occur in compact SL
$m$-folds in fairly generic almost Calabi--Yau $m$-folds.
A good test of this is whether the {\it deformation theory}
of compact SL $m$-folds with this kind of singularity is
well-behaved. That is, the analogue of Theorem \ref{gr5thm1}
should hold, with {\it finite-dimensional\/} obstruction
space~$\O_\sXp$.

One very obvious way to make examples of SL $m$-folds with
nonisolated singularities is to consider $C\t\R^{m-k}$ in
$\C^k\t\C^{m-k}=\C^m$, where $C$ is an SL cone in $\C^k$
with isolated singularity at 0, and $3\le k<m$. So we
could study SL $m$-folds with singularities locally
modelled on $C\t\R^{m-k}$. Calculations by the author
indicate that the deformation theory of such singular
SL $m$-folds will be well-behaved if and only if $C$ is
{\it stable}. Therefore we propose:

\begin{prob} Let $3\le k<m$, and suppose $C$ is an SL cone in
$\C^k$ with an isolated singularity at 0 which is {\it stable},
in the sense of Definition \ref{gr3def3}. Study compact SL
$m$-folds $N$ in almost Calabi--Yau $m$-folds $(M,J,\om,\Om)$,
where the singular set $S$ of $N$ is a compact
$(m\!-\!k)$-submanifold of $M$, and $N$ is modelled on
$C\t\R^{m-k}$ in $\C^k\t\C^{m-k}=\C^m$ at each singular
point~$s\in S$.
\label{gr8prob2}
\end{prob}

Here we have not defined what we mean by `modelled on'.
There should be some fairly natural asymptotic condition,
along the lines of \eq{gr3eq9}. Perhaps, as in Theorem
\ref{gr4thm3}, it will be equivalent to $N$ having tangent
cone $C\t\R^{m-k}$ with multiplicity 1 at each~$s\in S$.

A related problem is to classify the possible stable~$C$:

\begin{prob} Classify special Lagrangian cones $C$ in $\C^m$
for $m\ge 3$ with an isolated singularity at 0 which are
{\it stable}, in the sense of Definition~\ref{gr3def3}.
\label{gr8prob3}
\end{prob}

As above, by Haskins \cite[Th.~A]{Hask2} the SL $T^2$-cone
$C_{\sst\rm HL}^3$ in $\C^3$ of \eq{gr3eq4} is up to $\SU(3)$
equivalence the {\it unique} stable SL $T^2$-cone $C$ in $\C^3$.
In fact $C_{\sst\rm HL}^3$ is the {\it only} example of a
stable SL cone in $\C^m$ for $m\ge 3$ known to the author.
It is conceivable that it really is the only example, so that
the answer to Problem \ref{gr8prob3} is $C_{\sst\rm HL}^3$
and no others.

We can also look for other interesting classes of
singularities with well-behaved deformation theory.
The key is to find suitable asymptotic conditions.

\begin{prob} Let $C$ be an SL cone in $\C^m$ with nonisolated
singularity at 0, or with multiplicity $k>1$. Can you find a
good, natural set of asymptotic conditions for SL $m$-folds
with isolated singularities with tangent cone~$C$?
\label{gr8prob4}
\end{prob}

One way to approach this is through {\it examples}: we find
some class of examples of singular SL $m$-folds, calculate
their asymptotic behaviour near their singularities, and
try and abstract the important features. For the examples
in (i) above this may be easy, as they are very explicit.
But for those in (ii) above the author failed miserably to
understand the asymptotic behaviour.

\subsection{The SYZ Conjecture}
\label{gr83}

{\it Mirror Symmetry} is a mysterious relationship between pairs 
of Calabi--Yau 3-folds $M,\hat M$, arising from a branch of 
physics known as {\it String Theory}, and leading to some very 
strange and exciting conjectures about Calabi--Yau 3-folds.

Roughly speaking, String Theorists believe that each Calabi--Yau
3-fold $M$ has a quantization, a {\it Super Conformal Field
Theory} (SCFT). If $M,\hat M$ have SCFT's isomorphic under a
certain simple involution of SCFT structure, we say that
$M,\hat M$ are {\it mirror} Calabi--Yau 3-folds. One can argue
using String Theory that $H^{1,1}(M)\cong H^{2,1}(\hat M)$
and $H^{2,1}(M)\cong H^{1,1}(\hat M)$. The mirror transform
also exchanges things related to the complex structure of $M$
with things related to the symplectic structure of $\hat M$,
and vice versa.

The {\it SYZ Conjecture}, due to Strominger, Yau and Zaslow
\cite{SYZ} in 1996, gives a geometric explanation of Mirror
Symmetry. Here is an attempt to state it.

\begin{conj}[Strominger--Yau--Zaslow] Suppose $M$ and\/ $\hat M$ are
mirror Calabi--Yau $3$-folds. Then (under some additional conditions) 
there should exist a compact topological\/ $3$-manifold\/ $B$ and 
surjective, continuous maps $f:M\ra B$ and\/ $\hat f:\hat M\ra B$
with fibres $X_b=f^{-1}(b)$ and\/ $\hat X_b=\hat f^{-1}(b)$ for
$b\in B$, such that
\begin{itemize}
\setlength{\itemsep}{0pt}
\setlength{\parsep}{0pt}
\item[{\rm(i)}] There exists a dense open set\/ $B_0\subset B$, such 
that for each\/ $b\in B_0$, the fibres $X_b,\hat X_b$ are nonsingular
special Lagrangian $3$-tori $T^3$ in $M$ and\/ $\hat M$, which are in
some sense dual to one another.
\item[{\rm(ii)}] For each\/ $b\in\De=B\sm B_0$, the fibres $X_b$,
$\hat X_b$ are expected to be singular special Lagrangian $3$-folds
in $M$ and\/~$\hat M$.
\end{itemize}
\medskip
\end{conj}

We call $f,\hat f$ {\it special Lagrangian fibrations}, and the
set of singular fibres $\De$ is called the {\it discriminant}. It
is not yet clear what the final form of the SYZ Conjecture should
be. Much work has been done on it, working primarily with {\it
Lagrangian} fibrations, by authors such as Mark Gross and Wei-Dong
Ruan. For references see~\cite{Joyc5}.

The author's approach to the SYZ Conjecture, focussing primarily on
{\it special Lagrangian singularities}, is set out in \cite{Joyc5},
and we do not have space to discuss it here. Very briefly, we argue
that for generic (almost) Calabi--Yau 3-folds (ii) will not hold, as
the discriminants $\De,\hat\De$ of $f,\hat f$ cannot be homeomorphic
near certain kinds of singular fibre. We also suggest that the final
form of the SYZ Conjecture should be an {\it asymptotic statement\/}
about 1-{\it parameter families} of Calabi-Yau 3-folds approaching
the {\it large complex structure limit}.

\begin{prob} Study {\it special Lagrangian fibrations}
$f:M\ra B$ of almost Calabi--Yau 3-folds $(M,J,\om,\Om)$,
particularly when $\om$ is {\it generic} in its K\"ahler
class. Clarify/prove/disprove the SYZ Conjecture.
\label{gr8prob5}
\end{prob}

Note that the ideas of \S\ref{gr81} will be helpful here.
As $B$ has dimension 3, we see that $\ind(X_b)\le 3$ for
all $b\in\De$. If Conjecture \ref{gr5conj} holds, $\om$ is
generic, and $f^{-1}(b)$ has isolated conical singularities,
then $X_b$ is {\it transverse}. We can then use Theorem
\ref{gr8thm1} or \cite[Th.~8.10]{Joyc17} to calculate
$\ind(X_b)$, and $\ind(X_b)\le 3$ will severely restrict
the possible singular behaviour.

\subsection{Invariants from counting SL homology spheres}
\label{gr84}

In \cite{Joyc1} the author proposed to define an invariant of
almost Calabi--Yau 3-folds $(M,J,\om,\Om)$ by counting special
Lagrangian rational homology 3-spheres $N$ (which occur in
0-dimensional moduli spaces) in a given homology class, with
a certain topological weight. This invariant will only be
interesting if it is conserved under deformations of the
underlying almost Calabi--Yau 3-fold, or at least transforms
in a rigid way as the cohomology classes $[\om],[\Om]$ change.

During such a deformation, nonsingular SL 3-folds can develop
singularities and disappear, or new ones appear, which might
change the invariant. In \cite{Joyc1} the author showed that
if we count rational SL homology spheres $N$ with weight
$\bmd{H_1(N,\Z)}$, then under two kinds of singular behaviour
of SL 3-folds, the resulting invariant is independent of $[\om]$,
and transforms according to certain rules as $[\Om]$ crosses
real hypersurfaces in complex structure moduli space where
phases of $\al,\be\in H_3(M,\Z)$ become equal.
 
Again, the ideas of \S\ref{gr81} will be helpful here. It is
enough for us to study how the invariant changes along {\it
generic $1$-parameter families} of almost Calabi--Yau 3-folds.
The only kinds of singularities of SL homology 3-spheres that
arise in such families will have index 1. So if we can complete
the index 1 classification in Problem \ref{gr8prob1}, we should
be able to resolve the conjectures of~\cite{Joyc1}.

In fact, I now believe that interesting invariants of almost
Calabi--Yau $m$-folds by `counting' SL $m$-folds can be
defined for all $m\ge 3$. The definition, properties and
transformation laws of these invariants are formidably
complex and difficult, even to state. The best approach I
have to them is to use Homological Mirror Symmetry to translate
the problem to the derived category ${\cal T}=D^b($Fuk$(M,\om))$
of the Fukaya category of~$(M,\om)$.

Then SL $m$-folds conjecturally correspond to {\it stable
objects} of the triangulated category $\cal T$, under a
stability condition \`a la Tom Bridgeland. The invariants
are Euler characteristics of moduli spaces of {\it
configurations} in $\cal T$, which are finite collections
of (stable or semistable) objects and morphisms in $\cal T$
satisfying some axioms. In this set-up, using algebra and
category theory, I can rigorously develop the definition
and properties of the invariants, and their transformation
rules under change of stability condition (effectively,
deformation of $J,\Om$). I am writing (yet) another series
of papers about this.

\begin{prob} Try to use moduli spaces of compact SL
$m$-folds (possibly immersed, or singular) to define systems
of invariants of an almost Calabi--Yau $m$-fold $(M,J,\om,\Om)$
for $m\ge 3$. These invariants should be defined for $\om$ generic
in its K\"ahler class, and the key property we want is that they
should be {\it independent of\/} $\om$. Compute the invariants
for the quintic. Calculate the transformation rules for the
invariants under deformation of $J,\Om$. Relate them to
Homological Mirror Symmetry, and to `branes' in String Theory.
\label{gr5prob6}
\end{prob}

\bibliographystyle{amsplain}

\begin{thebibliography}{99}

\bibitem{Bred} G.E. Bredon, {\it Topology and Geometry}, Graduate
Texts in Mathematics 139, Springer-Verlag, Berlin, 1993.

\bibitem{Harv} R. Harvey, {\it Spinors and calibrations}, Academic Press, 
San Diego, 1990.

\bibitem{HaLa} R. Harvey and H.B. Lawson, {\it Calibrated geometries},
Acta Mathematica 148 (1982), 47--157.

\bibitem{Hask1} M. Haskins, {\it Special Lagrangian Cones},
math.DG/0005164, 2000. To appear in American Journal of Mathematics.

\bibitem{Hask2} M. Haskins, {\it The geometric complexity of
special Lagrangian $T^2$-cones}, math.DG/0307129, 2003. To
appear in Inventiones mathematicae.

\bibitem{Joyc1} D.D. Joyce, {\it On counting special Lagrangian homology
$3$-spheres}, pages 125--151 in {\it Topology and Geometry: Commemorating
SISTAG}, editors A.J. Berrick, M.C. Leung and X.W. Xu, Contemporary
Mathematics 314, A.M.S., Providence, RI, 2002. hep-th/9907013, 1999.

\bibitem{Joyc2} D.D. Joyce, {\it Special Lagrangian $m$-folds in $\C^m$ 
with symmetries}, Duke Math. J. 115 (2002), 1-51. math.DG/0008021.

\bibitem{Joyc3} D.D. Joyce, {\it Constructing special Lagrangian
$m$-folds in $\C^m$ by evolving quadrics}, Math. Ann. 320 (2001),
757--797. math.DG/0008155.

\bibitem{Joyc4} D.D. Joyce, {\it Evolution equations for special
Lagrangian $3$-folds in $\C^3$}, Ann. Global Anal. Geom. 20
(2001), 345--403. math.DG/0010036.

\bibitem{Joyc5} D.D. Joyce, {\it Singularities of special
Lagrangian fibrations and the SYZ Conjecture}, math.DG/0011179,
2000. To appear in Communications in Analysis and Geometry.

\bibitem{Joyc6} D.D. Joyce, {\it Ruled special Lagrangian $3$-folds
in $\C^3$}, Proceedings of the London Mathematical Society 85 (2002),
233--256. math.DG/0012060.

\bibitem{Joyc7} D.D. Joyce, {\it Special Lagrangian $3$-folds and
integrable systems}, 
\hfil\break
math.DG/0101249, 2001. To appear in volume 1 of the Proceedings of
the Mathematical Society of Japan's 9th International Research
Institute on {\it Integrable Systems in Differential Geometry},
Tokyo, 2000.

\bibitem{Joyc8} D.D. Joyce, {\it Lectures on Calabi--Yau and special
Lagrangian geometry},
math.DG/0108088, 2001. Published, with extra
material, as Part I of M. Gross, D. Huybrechts and D. Joyce,
{\it Calabi--Yau Manifolds and Related Geometries}, Universitext
series, Springer, Berlin, 2003.

\bibitem{Joyc9} D.D. Joyce, {\it $\U(1)$-invariant special Lagrangian
$3$-folds. I. Nonsingular solutions},
math.DG/0111324, 2001. To appear in Advances in Mathematics.

\bibitem{Joyc10} D.D. Joyce, {\it $\U(1)$-invariant special Lagrangian
$3$-folds. II. Existence of singular solutions}, math.DG/0111326, 2001.

\bibitem{Joyc11} D.D. Joyce, {\it $\U(1)$-invariant special Lagrangian
$3$-folds. III. Properties of singular solutions}, math.DG/0204343, 2002.

\bibitem{Joyc12} D.D. Joyce, {\it $\U(1)$-invariant special Lagrangian
$3$-folds in $\C^3$ and special Lagrangian fibrations}, Turkish
Math. J. 27 (2003), 99--114. math.DG/0206016.

\bibitem{Joyc13} D.D. Joyce, {\it Special Lagrangian submanifolds
with isolated conical singularities. I. Regularity}, math.DG/0211294, 2002.

\bibitem{Joyc14} D.D. Joyce, {\it Special Lagrangian submanifolds with
isolated conical singularities. II. Moduli spaces}, math.DG/0211295, 2002.

\bibitem{Joyc15} D.D. Joyce, {\it Special Lagrangian submanifolds with
isolated conical singularities. III. Desingularization, the unobstructed
case}, math.DG/0302355, 2003.

\bibitem{Joyc16} D.D. Joyce, {\it Special Lagrangian submanifolds with
isolated conical singularities. IV. Desingularization, obstructions
and families}, math.DG/0302356, 2003.

\bibitem{Joyc17} D.D. Joyce, {\it Special Lagrangian submanifolds with
isolated conical singularities. V. Survey and applications},
math.DG/0303272, 2003. To appear in the Journal of Diffierential Geometry.

\bibitem{Lawl} G. Lawlor, {\it The angle criterion}, Invent. math.
95 (1989), 437--446.

\bibitem{Lock} R. Lockhart, {\it Fredholm, Hodge and Liouville Theorems
on noncompact manifolds}, Trans. A.M.S. 301 (1987), 1--35.

\bibitem{Mars} S.P. Marshall, {\it Deformations of special Lagrangian
submanifolds}, Oxford D.Phil. thesis, 2002.

\bibitem{McSa} D. McDuff and D. Salamon, {\it Introduction to
symplectic topology}, second edition, OUP, Oxford, 1998.

\bibitem{McIn} I. McIntosh, {\it Special Lagrangian cones in $\C^3$
and primitive harmonic maps}, J. London Math. Soc. 67 (2003), 769--789.
math.DG/0201157.

\bibitem{McLe} R.C. McLean, {\it Deformations of calibrated submanifolds},
Communications in Analysis and Geometry 6 (1998), 705--747.

\bibitem{Morg} F. Morgan, {\it Geometric Measure Theory, A Beginner's
Guide}, Academic Press, San Diego, 1995.

\bibitem{Paci} T. Pacini, {\it Deformations of Asymptotically
Conical Special Lagrangian Submanifolds}, math.DG/0207144, 2002.

\bibitem{SYZ} A. Strominger, S.-T. Yau, and E. Zaslow, {\it Mirror 
symmetry is T-duality}, Nuclear Physics B479 (1996), 243--259.
hep-th/9606040.

\end{thebibliography}

\end{document}